\documentclass[a4paper]{article}
\usepackage{amsthm,amsmath,amssymb}
\usepackage[dvipsnames]{xcolor}
\newcommand{\fg}[1]{{#1}}

\newcommand{\ka}{{\kappa}}
\newcommand{\kc}{{k}}

\newcommand{\bele}{\begin{lemm}\begin{sl}}
\newcommand{\enle}{\end{sl}\end{lemm}}
\newcommand{\bedef}{\begin{defi}\begin{sl}}
\newcommand{\eddef}{\end{sl}\end{defi}}
\newcommand{\bete}{\begin{teor}\begin{sl}}
\newcommand{\ente}{\end{sl}\end{teor}}
\newcommand{\beos}{\begin{osse}\begin{rm}}
\newcommand{\eddos}{\end{rm}\end{osse}}
\newcommand{\bepr}{\begin{prop}\begin{sl}}
\newcommand{\empr}{\end{sl}\end{prop}}
\newcommand{\bepro}{\begin{prob}\begin{rm}}
\newcommand{\empro}{\end{rm}\end{prob}}
\newcommand{\bede}{\begin{defin}\begin{sl}}
\newcommand{\edde}{\end{sl}\end{defin}}
\newcommand{\beco}{\begin{coro}\begin{sl}}
\newcommand{\enco}{\end{sl}\end{coro}}

\newcommand{\quext}{\quad\text}

\newcommand{\RR}{\mathbb{R}}

\newcommand{\NN}{\mathbb{N}}

\newcommand{\beeq}[1]{\begin{equation}\label{#1}}
\newcommand{\eddeq}{\end{equation}}

\newcommand{\beeqa}[1]{\begin{eqnarray}\label{#1}}
\newcommand{\eddeqa}{\end{eqnarray}}

\newcommand{\beal}[1]{\begin{align}\label{#1}}
\newcommand{\eddal}{\end{align}}

\newcommand{\bespl}[1]{\begin{split}\label{#1}}
\newcommand{\edspl}{\end{split}}

\newcommand{\bega}[1]{\begin{gather}\label{#1}}
\newcommand{\edga}{\end{gather}}

\newcommand{\beeqax}{\begin{eqnarray*}}
\newcommand{\eddeqax}{\end{eqnarray*}}

\def\qed{\ifmmode % if math mode, assume display: omit penalty etc.
  \else \leavevmode\unskip\penalty9999 \hbox{}\nobreak\hfill
  \fi
  \quad\hbox{\hskip.5em\vrule width.4em height.6em depth.05em\hskip.1em}}
\def\endproofsym{\qed}
\def\endnobox{\def\endproofsym{}\end{proof}\def\endproofsym{\qed}}

\newcommand{\no}{\nonumber}

\newcommand{\unmezzo}{\frac{1}{2}}

\newcommand{\duav}[1]{\langle{#1}\rangle}

\newcommand{\perogni}{\forall\,}

\newcommand{\io}{\int_\Omega}

\newcommand{\epsi}{\varepsilon}

\newcommand{\OO}{_{\Omega}}

\newcommand{\lhs}{left-hand side}
\newcommand{\rhs}{right-hand side}

\DeclareMathOperator{\dive}{div}
\DeclareMathOperator{\deriv}{d}

\DeclareMathOperator{\sign}{sign}

\DeclareMathOperator{\reg}{reg}

\let\oldemph\emph
\renewcommand{\emph}[1]{\textbf{\oldemph{#1}}}

\let\TeXchi\chi
\def\chi{{\setbox0 \hbox{\mathsurround0pt
$\TeXchi$}\hbox{\raise\dp0 \copy0 }}}

\newcommand{\calD}{{\mathcal D}}

\newcommand{\calG}{{\mathcal G}}

\newcommand{\calE}{{\mathcal E}}

\newcommand{\dit}{\deriv\!t}
\newcommand{\dis}{\deriv\!s}
\newcommand{\dix}{\deriv\!x}

\newcommand{\ddt}{\frac{\deriv\!{}}{\dit}}
\newcommand{\derivt}{\frac{\deriv\!{}}{\dit}}

\newtheorem{thm}{Theorem}[]
\newtheorem*{thm-nn}{Theorem}

\newtheorem*{lem-nn}{Lemma}
\newtheorem{prop}[thm]{Proposition}

\newtheorem*{cor-nn}{Corollary}
\newtheorem{remark}{Remark}

\newcommand{\vt}{\vartheta}

\newenvironment{bettirev}{\color{blue}}{\color{black}}
\newcommand{\bber}{\begin{bettirev}}
\newcommand{\eber}{\end{bettirev}}

\newenvironment{michelarev}{\color{red}}{\color{black}}
\newcommand{\III}{\begin{michelarev}}
\newcommand{\EEE}{\end{michelarev}}
\newcommand{\vp}{\varphi}

\newcommand{\fhi}{\varphi}
\renewcommand{\u}{\textbf{u}}
\newcommand{\vu}{\textbf{u}}
\renewcommand{\v}{\textbf{v}}
\newcommand{\w}{\textbf{w}}
\newcommand{\grad}{\nabla} %gradiente 
\newcommand{\lap}{\Delta} %laplaciano
\renewcommand{\div}[1]{\nabla \cdot #1} %divergenza

\newcommand{\ede}[0]{:=}
\newcommand{\sprod}[2]{\langle #1 , #2 \rangle} %prodotto scalare
\newcommand{\weakconvstar}{\overset{*}{\rightharpoonup}} %convergenza debole-*
\newcommand{\weakconv}{\rightharpoonup} %convergenza debole
\newcommand{\immcont}{\hookrightarrow} %immersione continua
\newcommand{\dt}{\,\dit}
\newcommand{\dx}{\,\dix}

\DeclareMathOperator{\spann}{span}

\newcommand{\ints}[1]{\int_{\Omega}#1\dx} %integrale spaziale
\newcommand{\norm}[2]{\|#1\|_{#2}} %norma di #1 nello spazio #2
\renewcommand{\H}{H_{\dive}}
\newcommand{\V}{V_{\dive}}
\newcommand{\Hz}{H_{\dive,0}}
\newcommand{\Vz}{V_{\dive,0}}
\newcommand{\normHd}[1]{\|#1\|_{\H}} %norma di #1 nello spazio H_{\dive}
\newcommand{\normVd}[1]{\|#1\|_{\V}} %norma di #1 nello spazio V_{\dive}
\newcommand{\normH}[1]{\|#1\|_{H}} %norma di #1 nello spazio H
\newcommand{\normV}[1]{\|#1\|_{V}} %norma di #1 nello spazio V
\numberwithin{equation}{section}
\begin{document}

\title{On a Navier-Stokes-Allen-Cahn model with inertial effects}

\author{Gianluca Favre\\
Faculty of Mathematics, University of Vienna\\ 
Oskar-Morgenstern-Platz 1, 1090 Wien, Austria\\
E-mail: {\tt gianluca.favre@univie.ac.at}\\
\and
Giulio Schimperna\\
Dipartimento di Matematica, Universit\`a di Pavia,\\
Via Ferrata~5, I-27100 Pavia, Italy\\
E-mail: {\tt giusch04@unipv.it}
}

%\date{}

%\listofchanges
%\linenumbers

\maketitle
\begin{abstract}
 A mathematical model describing the flow of two-phase fluids in a 
 bounded container $\Omega$ is considered under
 the assumption that the phase transition is influenced by
 inertial effects. The model couples a variant
 of the Navier-Stokes system for the velocity $\vu$ with an
 Allen-Cahn-type equation for the order parameter $\fhi$ relaxed
 in time in order to introduce inertia at the microscopic level. 
 The resulting model
 is characterized by second-order material derivatives 
 which constitute the main difficulty in the mathematical analysis.
 In order to obtain a tractable problem, 
 a viscous relaxation term is then included in
 the phase equation. The mathematical results 
 consist in existence of weak solutions in 3D and, under additional
 assumptions, existence and uniqueness of strong solutions in 2D. 
 A partial characterization of the long-time behavior of solutions
 is also given and in particular some issues
 related to dissipation of energy are discussed.
\end{abstract}


\vspace{2mm}

\noindent {\bf AMS (MOS) subject clas\-si\-fi\-ca\-tion:}%
~~35Q35, 35K10, 35L82, 76D05, 80A22.

% % % % % % % % % % % % % % % % % % % % % % % % % % % % % % % % % % % % % % % % % % % % % % % % % % % % % %

\section{Introduction}

In this paper we are concerned with the mathematical analysis of the following
PDE system describing the evolution of a phase-changing fluid through the 
variables $\u$ (macroscopic velocity of the flow) and $\fhi$ (order parameter,
normalized in such a way that the values $\fhi=\pm1$ represent the pure states):
\begin{align}
\label{eq: N-S}
& \u_t + \u \cdot \grad \u + \grad p - \lap \u = - \dive(\grad\vp\otimes\grad\vp),\\
\label{eq: N-S div}
& \dive \u =0, \\
\label{eq: A-C pi}
& \pi = \vp_t + \u \cdot \grad\vp - \epsilon \lap\vp + \sigma \vp,\\
\label{eq: A-C pi_t}
& \ka\pi_t + \delta \u \cdot \grad\pi + \pi - \lap\vp + f(\vp) = 0,
\end{align}
with parameters $\ka > 0$, $\delta\ge 0$, $\sigma \ge 0$ and $\epsilon>0$ whose 
role will be clarified in the sequel. The evolution is assumed to take place in 
a bounded domain $\Omega\subset\RR^d$, $d\in\{2,3\}$, over a bounded but otherwise 
arbitrary time reference interval $(0,T)$. More precisely, in order to avoid complications
related with the interaction with the boundary, 
we will assume that $\Omega$ is the unitary flat torus 
$\Omega\sim [0,1]^d$. Correspondingly, we will take periodic boundary
conditions for all variables. This may be seen as the simplest situation
where one can attempt a mathematical analysis of the system; 
we however remark that at least a part of our results could be 
extended to other types of boundary conditions at the price of some additional 
technicality in the proofs. 

Relation \eqref{eq: N-S}, with the incompressibility constraint \eqref{eq: N-S div},
is the standard Navier-Stokes system where the 
\rhs\ of \eqref{eq: N-S} accounts for the extra-stress due to the effects
of phase transition. On the other hand, the change of phase
is governed by a proper version of the Allen-Cahn equation, where 
{\sl inertial effects}\/ are taken into account. Indeed, 
from a mechanical point of view, the phase transition phenomenon 
may be thought to be driven by the microscopic motion of molecules 
(we refer to the monograph \cite{Fre} for a general discussion),
which motivates the occurrence of an inertial term 
leading to a finite propagation speed of the phase transition front.
The function $f$ in \eqref{eq: A-C pi_t} represents the derivative of 
a configuration potential $F$ of double-well type whose minima are attained
in correspondence of the pure phases of the material. Mathematically
speaking, $F$ will be assumed to have a polynomial growth at infinity 
(of arbitrary order) and to be {\sl $\lambda$-convex}\/ (cf.~\eqref{powerlike}
below). It is worth remarking that, even
in dimension $d=2$, the polynomial behavior at infinity seems to be necessary 
in order for our arguments to work; in particular, our results do 
not extend to more general (e.g., logarithmic) potentials.

In the above formulation, the Allen-Cahn equation for the phase parameter has been split
into the coupled relations \eqref{eq: A-C pi}, \eqref{eq: A-C pi_t}.
Indeed, this choice seems more natural because it gives a better description
of the transport effects in terms of the auxiliary variable $\pi$,
which plays an important role from both the physical and the 
mathematical viewpoint. Indeed, from \eqref{eq: A-C pi} we can see 
that $\pi$ represents a sort of regularization of the material
derivative of $\fhi$, $\epsilon$ and $\sigma$ being the (small)
regularization parameters. The reasons for considering such regularization
terms are shortly discussed in the next paragraphs, whereas a much more detailed 
explanation (also considering the physical implications) is given in
Section~\ref{sec:conrem} below. Here, it is also worth observing that,
for $\ka = \delta > 0$, equation \eqref{eq: A-C pi_t} 
includes the term $\pi_t + \u \cdot \grad \pi$,
which, by \eqref{eq: A-C pi}, can be interpreted as a regularization of 
the second {\sl material}\/ derivative of $\fhi$.

The mathematical literature devoted to the study of models for complex 
(i.e., two-phase or two-component) fluids is very vast. Here we limit ourselves
to describe a limited number of contributions without any claim of completeness.
The bibliographic entries in the papers mentioned below can provide the reader
with a further set of references. The first article presenting a mathematical 
model for two-phase fluids, dating back to 1977, is \cite{Hohenberg-Halperin'77},
where the so-called ``H-model'' for binary fluids was introduced. Actually,
the approach devised there is still at the base of more modern models.
Later, in 1996, Gurtin et al.~\cite{Gurtin et al.'96} introduced a coupled Navier-Stokes-Cahn-Hilliard 
whose physical derivation was obtained using a balance law for microforces in conjunction with 
constitutive equations consistent with a mechanical version of the second law of Thermodynamics.
This mechanical interpretation of the phase transition effect as an action
of microscopic forces is also at the base of the present model and, in particular,
justifies inertial effects at the microscopic level. A comprehensive presentation
of mechanically-based phase change models is given in the monograph
\cite{Fre}. A situation where the (compressible) Navier-Stokes 
system is coupled with an equation of Allen-Cahn type
was then introduced in \cite{Bles} in order to provide a mathematical description of 
the physical problem of cavitation in a flowing fluid. In particular, in that
contribution, a Korteweg term having the same expression $-\dive(\grad\vp \otimes \grad\vp)$ as in
our~\eqref{eq: N-S} appears as a forcing term on the \rhs\ of the Navier-Stokes system.

The above mentioned papers are mainly devoted to presenting a number of significant
models and discussing their compatibility with the basic laws of Thermodynamics and Fluid Mechanics.
Speaking, instead, of mathematical results, the literature devoted to the analysis 
(also from the numerical viewpoint) of models for complex fluids is more recent.
In 2003, Liu and Shen \cite{LS} provided a variational formulation of a 
Cahn-Hilliard-Navier-Stokes for binary fluids with matched densities and viscosities
and proved analytical and numerical results partly based on a Fourier-spectral method.
In 2009, Abels \cite{Ab} studied a similar
model coupling the incompressible Navier-Stokes system with the Cahn-Hilliard equation 
and proved existence both of weak and of strong solutions under natural conditions on data. 
In particular he considered the case when the two fluids have the same density but
may have different viscosities. He also discussed the sharp-interface limit of the resulting model.
A further related contribution \cite{AbelsRoger'09} is devoted to the analysis of a 
Navier-Stokes-Mullins-Sekerka model. Another important class of results is devoted 
to studying the long-time properties of solutions; 
among the many related works, we may mention \cite{Zhao Wu Huang'09}, 
where the  long time behavior of the strong solution to a model of viscous incompressible fluids
is analyzed. In particular the rate of the convergence to equilibria is established there. 
More recently, in \cite{GG10}, the existence of global and exponential
attractors has been established for a Cahn-Hilliard-Navier-Stokes model
in two space dimensions. Finally, we may mention \cite{CG12}, where the global regularity 
of strong solutions to a Cahn-Hilliard-Navier-Stokes
model with mixed partial viscosity and mobility is analyzed. 

Second order in time (hyperbolic or pseudo-parabolic)
relaxations of the Allen-Cahn and Cahn-Hilliard
equations have been considered in several different contexts
and under various mathematical assumptions. We may mention,
without any claim of completeness, the papers \cite{Bonfanti2004123,GP1,GMPZ,GSZ2D}
(see also the references therein). In most of these contributions, 
thermal effects have also been considered, leading to the so-called {\sl phase-field}\/
models with inertia. Inertial effects at the microscopic 
level have been observed also in (nematic)
liquid crystals (see, e.g., \cite{FRSZhyp,Gay-Balmaz}). On the other
hand, at least up to our knowledge, phase transition models with inertia
have not been studied so far in the case of two-phase
fluids, i.e., when the process is influenced by a macroscopic
velocity satisfying some version of the Navier-Stokes system.

Coming to a description of our results, we first observe that
the term $- \epsilon\lap\vp$ in \eqref{eq: A-C pi}, 
acts as a viscous regularization of the phase variable. 
This type of effect may be experimentally observed
at least as far as $\epsilon>0$ is very small, and 
may correspond to a {\sl strong damping}\/ in the phase transition.
On the other hand, here such term is added mainly for
mathematical reasons (i.e., it acts as a smoothing
parameter). Indeed, omitting it (i.e., taking $\epsilon=0$),
the resulting system seems to present, even in 2D, insurmountable
mathematical difficulties, especially related to the control 
of the \rhs\ of \eqref{eq: N-S} (see also Remark~\ref{rem:eps} 
below for additional comments).

Our existence proof will be based on a standard Faedo-Galerkin
regularization scheme, combined with the derivation of 
a priori estimates and compactness arguments to pass to the limit.
Entering details, we start observing that the only 
{\sl a priori}\/ information that is always guaranteed for solutions to our system
is the {\sl energy estimate}, which is a direct consequence of 
the variational structure of the model (cf.~Subsec.~\ref{subsec:en} 
below for more details). 
On the other hand, this bound is sufficient only
for proving existence of weak solutions. In particular, it is only
thanks to the viscous regularization in \eqref{eq: A-C pi} that
we can obtain some information on the second space derivatives of $\fhi$.
Actually, for $\epsilon=0$ taking the limit of
the right-hand side of \eqref{eq: N-S} would likely be hopeless 
due to lack of compactness. Indeed, here
it does not seem possible to obtain the strong $L^2$-convergence of $\nabla\fhi$ 
via contractive arguments (as is usually done for the decoupled semilinear wave
equation with subcritical nonlinearity), due to the presence of the 
transport terms.

Once we have obtained existence of weak solutions (which we 
can do both in 2D and in 3D) 
looking for additional properties is more difficult, even in the 2D case,
and we can only prove a number of partial results 
holding under additional regularity assumptions
on coefficients and data. Entering details, in 2D we can try to construct 
{\sl strong solutions}\/ (such a classification is mutuated 
by the standard terminology used for Navier-Stokes)
and, indeed, we can prove their existence, but only in the case when the second
order transport term is neglected, i.e., one has $\delta=0$ in \eqref{eq: A-C pi_t}. 
In the class of strong solutions one can also prove uniqueness. 
This result holds, with a conditional nature, also in 3D
(where existence of strong solutions is not known, of course),
or for $\delta>0$. 

A further question we address is related to energy dissipation. Actually,
due to the absence of external sources and to the choice of 
periodic boundary conditions,
one expects that at least a part of the kinetic and chemical energy
of the body is gradually converted into heat. On the other hand, we
can explicitly prove this fact only in two cases, i.e.~when either
$\sigma>0$ in \eqref{eq: A-C pi} or $\delta = \ka$ in \eqref{eq: A-C pi_t}. 
Actually, if $\sigma>0$ one can control, uniformly in time, 
the (possibly nonpositive) contribution coming from the combination 
of the additional viscosity term $-\epsilon\Delta \fhi$
with the {\sl non-monotone}\/ semilinear term $f(\fhi)$.
On the other hand, for $\delta = \ka$, corresponding to the 
situation when the true material derivative of $\pi$ occurs 
in \eqref{eq: A-C pi_t}, dissipativity holds
even for $\sigma = 0$ thanks to a cancellation property.
Under different choices on regularization coefficient, 
dissipativity of energy may actually not hold.
This (somehow unexpected) behavior is correlated with the presence 
of the viscous regularization parameter $\epsilon$ 
in \eqref{eq: A-C pi}. A discussion
of this issue, and of its physical implications, is given in 
Section~\ref{sec:conrem} below. There, we also mention the possibility
of considering, in place of $\sigma\fhi$, a more general {\sl semilinear}\/ 
contribution $\sigma(\fhi)$ in \eqref{eq: A-C pi} (with one relevant 
choice being given by $\sigma(\cdot)=f(\cdot)$).
This term, at least under suitable conditions on the 
function $\sigma$, would describe some form of 
{\sl nonlinear damping}\/ effect (see, e.g., \cite{WD-5 Barbu, WD-4 Bociu, 
WD-2Chueshov2004, WD-7 Dell'Oro, WD-6 Dell'Oro, Khanmamedov, WD1-LASIECKA, Roder-Tebou} 
for related models); the resulting mathematical problem
will be possibly addressed in a forthcoming paper.

It is finally worth noting that, in the language of dynamical systems,
uniform dissipativity of energy implies the 
existence of a {\sl bounded absorbing set}\/ in a suitable phase-space,
a fact that may serve as a starting point for studying the long-time behavior
of solution trajectories and proving existence of attractors.
These issues will be also possibly addressed in a future work.

\smallskip

We conclude with the plan of the paper. The next section is devoted to 
presenting our precise mathematical assumptions on coefficients and data
and a number of preliminary considerations. Then, our main results
are stated in the subsequent Section~\ref{sec:main}. The core
of the proofs, including the basic a priori estimates and the 
arguments used for passing to the limit in the approximation, are
then presented in Section~\ref{sec:pro}. A possible construction
of regularized solutions by means of a Faedo-Galerkin scheme
is given in Section~\ref{sec:Gal app}. Finally, in Section~\ref{sec:conrem},
we provide some final comments and, in particular, we discuss
a bit more extensively the physical implications of our choices
about regularizing parameters.

% % % % % % % % % % % % % % % % % % % % % % % % % % % % % % % % % % % % % % % % % % % % % % %

\section{Preliminaries}

% % % % % % % % % % % % % % % % % % % % % % % % % % % % % % % % % % % % % % % % % % % % % % %

\subsection{Notation and functional setup}

% % % % % % % % % % % % % % % % % % % % % % % % % % % % % % % % % % % % % % % % % % % % % % %

We will note as $\Omega$ the unit flat torus $[0,1]^d$, with $d=2$ or $d=3$. As is customary,
all functions defined on $\Omega$ will be implicitly assumed to satisfy periodic boundary
conditions in a suitable sense. This will not be emphasized in the notation, for 
the sake of simplicity. For instance, we will set $H:=L^2(\Omega)$ and $V:=H^1(\Omega)$ 
implicitly assuming $\Omega$-periodicity. These spaces will be used as function
spaces for the phase variable; the same symbols $H$ and $V$
will be used also for denoting vector- or tensor-valued functions (we may write, for instance,
$\nabla\fhi \in H$). The standard scalar product in $H$ will be noted as 
$(\cdot,\cdot)$. Since the immersion $V\subset H$ is continuous
and dense, identifying $H$ with $H'$ through the above scalar product
we obtain the {\sl Hilbert triplet}\/ $(V,H,V')$ for the phase variable.

Concerning the velocity function $\u$, we set
\begin{equation}\label{spazi:u}
  \mathbf{C}^\infty_{\dive}(\Omega) = \big\{\u \in \big[C^\infty(\Omega)\big]^d : \dive \u = 0 \big\}, \qquad
    \V:= \overline{\mathbf{C}^\infty_{\dive}(\Omega)}^{H^1(\Omega)}, \quad
    \H := \overline{\mathbf{C}^\infty_{\dive}(\Omega)}^{L^2(\Omega)},
\end{equation}
where, again, $\Omega$-periodicity is still implicitly subsumed everywhere.
The spaces $\H$ and $\V$ are seen as (closed) subspaces of $H$ and $V$ (more precisely, of $H^d$ 
and $V^d$), respectively, and in particular they are endowed with the corresponding norms.
Then, the embedding $\V\subset \H$ is continuous and dense, which permits us to identify $\V$
with its topological dual $\V'$ by means of the scalar product of $\H$, still denoted by
$(\cdot,\cdot)$, and to construct the ``velocity Hilbert triplet'' $\big(\V,\H,\V'\big)$.
We recall that $\V$ is endowed with the bilinear form $(\!(\u,\v)\!) \ede (\grad\u,\grad\v)$
for all $\u,\v \in \V$. Moreover, given a generic Banach space $X$ (and particularly in the
cases $X=V$, $X=\V$), we will denote by $\duav{\cdot,\cdot}$ the duality
between $X'$ and $X$. 
As is customary, we define the trilinear $\V$-continuous form
\begin{equation}
\label{def: trilinear b}
b(\u,\v,\w) = \ints{(\u\cdot\grad \v)\cdot\w} \quad \perogni \u,\v,\w \in \V
\end{equation}
and the bilinear form $\mathbf{B}:\V\times\V\to\V'$ given by
\begin{equation*}
\sprod{\mathbf{B}(\u,\v)}{\w}=b(\u,\v,\w) \quad \perogni \u,\v,\w \in \V.
\end{equation*}
Using Ladyzhenskaya's inequality (see, e.g. \cite[Chap.~5]{Robinson}), 
it turns out that the following properties hold for every $\u$, $\v$ and $\w$~$\in\V$:
%
%\begin{gather}
%\label{eq: b lim1}
%|b(\u,\v,\w)| \le c \norm{\u}{\H}^{1/2}\norm{\grad\u}{\H}^{1/2}\norm{\grad\v}{\H}\norm{\grad\w}{\H} \quad \dim(\Omega)=3\\
%\label{eq: b lim2}
%|b(\u,\v,\w)| \le c \norm{\u}{\H}^{1/2}\norm{\grad\u}{\H}^{1/2}\norm{\grad\v}{\H}\norm{\w}{\H}^{1/2}\norm{\grad\w}{\H}^{1/2} \quad \dim(\Omega)=2
%\end{gather}
%
\begin{gather}
 \label{eq: b lim1}
  |b(\u,\v,\w)| \le c \norm{\u}{\H}^{1/2}\norm{\u}{\V}^{1/2}\norm{\v}{\V}\norm{\w}{\V} \quext{if }\,d=3,\\
 \label{eq: b lim2}
  |b(\u,\v,\w)| \le c \norm{\u}{\H}^{1/2}\norm{\u}{\V}^{1/2}\norm{\v}{\V}\norm{\w}{\H}^{1/2}\norm{\w}{\V}^{1/2} \quext{if }\,d=2.
\end{gather}

% % % % % % % % % % % % % % % % % % % % % % % % % % % % % % % % % % % % % % % % % % % % % % %

\subsection{Assumptions on the potential}

We let $F\in C^2(\mathbb{R};\mathbb{R})$ and setting $f \ede F'$ we assume that
\begin{equation}\label{powerlike}%\tag{A1}
  \fg{K_0}(|s|^p + 1) \ge f'(s) \ge k_0|s|^p - \fg{\lambda_0}
\end{equation}
for some $\fg{K_0, k_0} , p >0$, $\fg{\lambda_0} \ge 0$ and every $s\in \mathbb{R}$. Namely, $f$ grows at infinity 
as a power-like function of exponent $p+1>1$ (i.e., it is strictly superlinear). Note that
the above implies in particular that $F$ is (at least) $\lambda$-convex, i.e., $F''=f'$ is
everywhere greater than $\fg{-\lambda_0}$, with true convexity holding in the case $\fg{\lambda_0 = 0}$.

A practical and common example of a phase potential satisfying the above is the 
standard \emph{double well potential}\ having the expression
$F(s)= (s^2 - 1)^2$ for which \eqref{powerlike} holds of course for $p=2$.

The above relation has several useful consequences, some of which are listed below. 
First, by integration one can easily deduce 
\begin{equation}\label{powerlike2}
  \fg{K'}(|s|^{p+1} + 1) \ge f(s)\sign s \ge \fg{k'}|s|^{p+1}  - \fg{\lambda_1},
\end{equation}
and, integrating again, 
\begin{equation}\label{powerlike3}
  \fg{K''}(|s|^{p+2} + 1) \ge F(s) \ge \fg{k''} |s|^{p+2}  - \fg{\lambda_2},
\end{equation}
where $\fg{K', k', K'', k'' >0}$ and $\fg{\lambda_1,\lambda_2 \ge 0}$ are computable constants
depending only on $\fg{K_0, k_0 , p , \lambda_0}$. Moreover, combining 
\eqref{powerlike2} and \eqref{powerlike3}, it is not difficult to deduce that
\begin{equation}\label{frr}
  \kc_p ( |s|^{p+2} + F(s) ) - c_p
   \le f(s) s \le C_p ( |s|^{p+2} + 1 ) \quad \perogni s \in \RR,
\end{equation}
where the constants $\kc_p, C_p >0$ and $c_p\ge 0$ depend only on the given values of the
parameters in~\eqref{powerlike}. Moreover, thanks to $p+2>2$, we also observe
that for any (large) $M>0$ and (small) $\epsi>0$ there exists $c(M,\epsi)>0$
such that
\begin{equation}\label{frr2}
  M | s |^2 \le \epsi F(s) + c(M,\epsi)  \quad \perogni s \in \RR.
\end{equation}

% % % % % % % % % % % % % % % % % % % % % % % % % % % % % % % % % % % % % % % % % % % % % % %

% % % % % % % % % % % % % % % % % % % % % % % % % % % % % % % % % % % % % % % % % % % % % % %

\subsection{Initial data and conservation properties}

In the sequel we will note the spatial mean of a generic
function $v$ defined over $\Omega$ as
\begin{equation}\label{def:media}
  v\OO := \frac{1}{|\Omega|} \io v \,\dix = \io v \,\dix,
\end{equation}
the second equality holding because $\Omega$ is the {\sl unit}\/ torus.
We also recall the Poincar\'e-Wirtinger inequality 
\begin{equation}\label{powi}
  \| v - v\OO \|_H \le c \| \nabla v \|_H,
\end{equation}
holding for any $v \in V$ and for some $c>0$. In particular, 
thanks to periodic boundary conditions, the
above holds with $v\OO = 0$ when $v = D_{x_i} z$ for
some $z\in H^2(\Omega)$ and $i\in\{1,\dots,d\}$.

We observe that, due to the choice of periodic boundary conditions, any hypothetical solution 
to our problem satisfies some conservation properties. First of all,
integrating \eqref{eq: N-S} over $\Omega$, we actually have
\begin{equation}\label{cons:u}
  \ddt \u\OO = 0,
\end{equation}
i.e., the mean value of the velocity is conserved in time. Such a property corresponds to
the conservation of (total) momentum in absence of external forces. 

Then, in order to reduce technical complications
we shall always assume that
\begin{equation}\label{hp:u}
  \u_0 \in \H, \quad (\u_0)\OO = 0.
\end{equation}
Indeed the case when $(\u_0)\OO \neq 0$ could be reduced to the present one by 
rewriting the system in terms of the translated variable $\u - \u\OO$. 

Next, integrating \eqref{eq: A-C pi} and \eqref{eq: A-C pi_t}, we deduce, respectively,
\begin{align}\label{cons:fhi1}
  & \ddt \fhi\OO + \sigma \fhi\OO = \pi\OO, \\
 \label{cons:fhi2}
  & \ka \ddt\pi\OO + \pi\OO + \io f(\fhi) \,\dix = 0,
\end{align}

whence, combining the above relations, we obtain (here primes denote
derivation in time)
\begin{equation}\label{cons:fhi}
  \ka  (\fhi\OO)'' + (1+\ka \sigma)(\fhi\OO)' + \sigma \fhi\OO + \io f(\fhi) \,\dix = 0.
\end{equation}
This relation rules the evolution of the total mass of either
component of the binary fluid. Note, however, that a decay property
for $\fhi\OO$ does not follow directly from \eqref{cons:fhi}
due to the presence of the nonlinear function $f$.

%%%%%%%%%%%%%%%%%%%%%%%%%%%%%%%%%%%%%%%%%%%%%%%%%%%%%%%%%%%%%%%%%%%%%%%%%%%%%%%%%%%%%%%%%%%%%%%%%%%%%%%%%%%%%%%%%%%%%%%%%%%%%%%%%%%%%%

\section{Main Results}
\label{sec:main}

Basically we can distinguish our results into two classes. The first one refers to the regularity setting
of {\sl weak solutions}\/ (this essentially means that the initial data have the sole regularity
corresponding to the finiteness of the physical energy). In such a framework, we can prove 
existence both in 2D and in 3D (Theorem~\ref{thm: esistenza debole}) for $\delta\ge 0$. Moreover, 
in all cases except when $\sigma = 0$ and $\delta \ne \ka$, we can also prove an
energy dissipation principle (Prop.~\ref{thm: dissipativita}), i.e.~the 
fact that, whatever is the magnitude of the energy at the initial time, solution trajectories
(or, at least those solution trajectories that arise as cluster points of approximating families;
indeed, uniqueness is not known at this level) tend to have an energy configuration that is below
some computable threshold depending only on the parameters of the system (and not on the 
magnitude of the initial data). In the language of dynamical systems this corresponds
to the existence of a {\sl uniformly absorbing set}. The reasons why we are able to
prove such a condition only under additional assumptions on coefficients are detailed
in Section~\ref{sec:conrem} below.

The second family of results refer to the class of {\sl strong solutions}, hence requiring
some more smoothness of initial data. Unfortunately, even in 2D, we can prove existence of such solutions
(Theorem~\ref{thm: esistenza forte}) only in the case when $\delta=0$, i.e.~the second-order transport
effect is neglected. In this case, we can also prove uniqueness of strong solutions 
(Theorem~\ref{thm: unicit}). The result holds both in 2D and in 3D, having of course a conditional 
nature in the latter case.
%
%\begin{remark}
% There is an interesting and unusual point in these approach because we can observe how transport terms play a fundamental role 
% to define the dissipativity of the system. A priori it is not expected that the dissipation can depend from transport terms.
% If we consider both the terms we will show that doing the dissipativity estimate they erase each others. Otherwise, studying 
% the \emph{simplified} system, the transport term $\u\cdot\grad\vp$ does not allow us to find an in time compact set for $\normH{\vp(t)}$. 
%\end{remark}
%
%\begin{remark}\label{rem:dissi}{\rm
% %
% When $\sigma = 0$, the lack of a dissipative estimate in the 
% general case is tied to the occurrence of the viscosity term $-\epsilon \lap\vp$ 
% in \eqref{eq: A-C pi}. Indeed, at least in the physical case when $F$ is nonconvex, 
% the energy inequality (cf.~\eqref{eq:energy} below) may contain a nonpositive term on
% the \lhs\ and also the information
% \eqref{cons:fhi} on the spatial means seems not sufficient to control it uniformly in time.
% Hence an estimate of the energy can be obtained only using Gronwall's inequality, which
% may lead to an exponential growth in time. On the other hand, in order to avoid this (unexpected) 
% behavior, we may try to adapt the argument standardly used to prove dissipativity of solutions to the 
% semilinear wave equation. However, we then face the occurrence of the transport terms, which
% apparently cannot be controlled in the general case, whereas they
% vanish thanks to incompressibility if $\delta = \ka$ (and only in that case).}
 %
%\end{remark}
%
%\noindent%

That said, we start with stating our first result devoted to existence of weak solutions:
\begin{thm}[Existence of weak solutions, $d=2,3$] \label{thm: esistenza debole}
 Let $d\in\{2,3\}$, let $\epsilon>0$, $\sigma \ge 0$, $\ka > 0$ and $\delta\ge 0$. Let also Assumption~\eqref{powerlike}
 hold. Let $T>0$ and let the initial data satisfy
 \begin{equation}\label{hp:init}
   \u_0\in\H, \quad (\u_0)\OO = 0, \qquad 
    \vp_0\in V \cap L^{p+2}(\Omega), \qquad
    \pi_0\in H.
 \end{equation}
 Then there exists at least one triple $(\u,\vp,\pi)$ belonging to the 
 regularity class
 \begin{align}\label{rego:u:w}
   & \u \in W^{1,q}(0,T;\V')\cap L^\infty(0,T;\H) \cap L^2(0,T;\V),\\
  \label{rego:fhi:w}  
   & \vp \in L^\infty(0,T;V)\cap L^\infty(0,T;L^{p+2}(\Omega)) \cap L^2(0,T;H^2(\Omega)),\\
  \label{rego: fhi t}
   & \vp_t \in L^2(0,T;L^s(\Omega)),\\
  \label{rego:pi:w} 
   & \pi \in L^\infty(0,T;H),\\
  \label{rego: pi t}
   & \pi_t \in L^2(0,T;(W^{1,3})'(\Omega)),
 \end{align}
 where 
 \begin{equation}\label{espon}
    q=s=2 ~~\text{if }\,d=2, \qquad 
    q=4/3,~ s=3/2 ~~\text{if }\,d=3,
 \end{equation}
 satisfying system~\eqref{eq: N-S}-\eqref{eq: A-C pi_t} in the following weak sense
 \begin{align}\label{N-S:weak} 
   & \duav{ \u_t , \v } 
    - \io (\u\otimes\u) : \nabla \v \, \dix
    + ( \nabla \u, \nabla \v )
   = \io ( \nabla \fhi \otimes \nabla \fhi ) : \nabla \v \, \dix,\\
  \label{AC:pi:weak} 
   & \pi = \vp_t + \u \cdot \grad\vp - \epsilon\lap\vp + \sigma \vp, 
    \quad\text{a.e.~in }\Omega,\\
  \label{AC:pi_t:weak} 
   & \ka \duav{\pi_t, v } 
    + ( \pi, v ) 
    - \delta \io \pi \u \cdot \nabla v \, \dix
    + ( \nabla \fhi, \nabla v )
    + \io f(\fhi) v \, \dix
    = 0,
 \end{align}
 for almost every $t\in (0,T)$, every $\v \in \V$ and every $v \in W^{1,3}(\Omega)$,
 and complying with the initial conditions
 \begin{equation}\label{init}
   \u|_{t=0} = \u_0, \qquad
    \fhi|_{t=0} = \fhi_0, \qquad
    \pi|_{t=0} = \pi_0    
 \end{equation}
 almost everywhere in $\Omega$. The triple $(\u,\vp,\pi)$ will be noted as a\/ {\rm weak
 solution} in the sequel.
\end{thm}
\begin{remark}\label{rem:en}{\rm
 It is worth observing that condition \eqref{hp:init} corresponds exactly to the finiteness
 of the physical energy (cf.~\eqref{eq:energy} below) 
 at the initial time. Indeed, assumption $\fhi_0\in L^{p+2}(\Omega)$ is equivalent
 to asking $F(\fhi_0)\in L^1(\Omega)$ due to \eqref{powerlike}. }
\end{remark}
\noindent%
Since $T>0$ is arbitrary, we can assume that the weak solutions
provided by Theorem~\ref{thm: esistenza debole} could be extended to be defined
for any time $t\in[0,\infty)$. In this perspective we can prove the following
\begin{prop}[Dissipativity] \label{thm: dissipativita}
 Let the hypotheses of Theorem~\ref{thm: esistenza debole} hold and let us 
 additionally assume that, either $\sigma > 0$, or $\delta = \ka$.
 Then there exists a constant $C_0$ independent of the initial data
 (but depending on the other parameters of the system) and a 
 time $T_0$ depending only on the ``energy'' of the initial datum,
 i.e.~on the quantity
 \begin{equation}\label{magn:init}
   \| \u_0 \|_{\H}
    + \| \fhi_0 \|_{V}
    + \| \fhi_0 \|_{L^{p+2}(\Omega)}
    + \| \pi_0 \|_{H},
 \end{equation}
 such that any weak solution
 provided by Theorem~\ref{thm: esistenza debole} emanating from this initial 
 datum satisfies
 \begin{equation}\label{st:diss}
   \| \u(t) \|_{\H}
    + \| \fhi(t) \|_{V}
    + \| \fhi(t) \|_{L^{p+2}(\Omega)}
    + \| \pi(t) \|_{H}
   \le C_0 \quad \perogni t \ge T_0.
 \end{equation}
 %
 %
 %In particular a good choice for $T_0$, considering $C_0(\sigma, c_p, \gamma)$ as in 
 %\eqref{insieme assorbente}, is $T_0 = \max \bigg\{0, \frac{\ln \calG(0)}{\gamma} \bigg\}$.  
 %
\end{prop}
\begin{remark}\label{T:grande}{\rm
 In the terminology of dynamical systems, the constant $C_0$ 
 in estimate \eqref{st:diss} may be interpreted 
 as the ``radius'' of an absorbing set with respect to the norms specified 
 there (which in turn somehow quantify the magnitude of the energy). }
\end{remark}
\begin{remark}\label{on:uniq} {\rm 
 Since the uniqueness of weak solutions is not known at this level, the dissipativity property
 in \eqref{st:diss} has to be carefully interpreted. Indeed, its proof is obtained by passing
 to the limit in an analogue relation holding at the approximate level (i.e., for some regularized
 solution that has better properties). Hence, any weak solution that is a limit point of a 
 sequence of approximate solutions turns out to satisfy it. On the other hand, we cannot
 exclude that there might exist ``bad'' weak solutions, unrelated to the approximation
 scheme, that do not satisfy \eqref{st:diss}. }
\end{remark}
\begin{thm}[Existence of strong solutions]
 \label{thm: esistenza forte}
 Let $d=2$, $\epsilon > 0$, $\kappa >0$, $\delta = 0$, $\sigma \ge 0$ and let Assumption~\eqref{powerlike}
 hold. Let also $T>0$ and let the initial data satisfy
 \begin{equation}\label{hp:init2}
   \u_0\in\V, \qquad 
   (\u_0)\OO = 0, \qquad 
   \vp_0\in H^2(\Omega), \qquad
    \pi_0\in V.
 \end{equation}
 Then there exists almost one {\sl strong solution}, namely a triple $(\u,\vp,\pi)$
 with 
 \begin{align}\label{rego:u:s}
   &  \u \in H^1(0,T;\H)\cap L^\infty(0,T;\V) \cap L^2(0,T;H^2_{\dive}(\Omega)),\\
  \label{rego:fhi:s}  
   &  \vp \in H^1(0,T;V) \cap L^\infty(0,T;H^2(\Omega))\cap L^2(0,T;H^3(\Omega)),\\
  \label{rego:pi:s}
   &  \pi \in H^1(0,T;H) \cap L^\infty(0,T;V),
 \end{align}
 satisfying system\/~\eqref{eq: N-S}-\eqref{eq: A-C pi_t} in the following sense:
 \begin{align}\label{NS-strong}
  & \u_t + \dive(\u\otimes\u) + \grad p + \dive(\grad\vp\otimes\grad\vp) - \lap \u = 0
   \quext{a.e.~in }\,(0,T)\times\Omega,\\
 \label{div-strong}
  & \dive \u =0\quext{a.e.~in }\,(0,T)\times\Omega,\\
 \label{pi-strong}
  & \pi = \vp_t + \u \cdot \grad\vp - \epsilon\lap\vp + \sigma \vp 
  \quext{a.e.~in }\,(0,T)\times\Omega,\\
 \label{fhi-strong}
  & \kappa \pi_t + \pi - \lap\vp + f(\vp) = 0\quext{a.e.~in }\,(0,T)\times\Omega
 \end{align}
 and the initial conditions as in \eqref{init}.
\end{thm}
\begin{thm}[Uniqueness of strong solutions]
 \label{thm: unicit}
 Let $d\in\{2,3\}$, $\sigma \ge 0$, $\kappa >0$,
 $\epsilon > 0$ and let Assumption~\eqref{powerlike} hold.
 Let also $T>0$ and let $(\u_i,\vp_i,\pi_i)$, $i=1,2$ be a couple of
 strong solutions (in the sense of the previous theorem, and satisfying
 in particular the regularity conditions\/ \eqref{rego:u:s}-\eqref{rego:pi:s})
 emanating from the same initial datum $(\u_0,\vp_0,\pi_0)$ 
 satisfying \eqref{hp:init2}. Moreover, let either $\delta = 0$ 
 or $\delta>0$ together with the additional assumption
 \begin{equation}\label{rego:add}
   \pi_i\in L^2(0,T;W^{1,3}(\Omega)) \quext{for some }\,i\in\{1,2\}.
 \end{equation}
 Then $(\u_1,\vp_1,\pi_1)=(\u_2,\vp_2,\pi_2)$ almost everywhere in 
 $(0,T)\times \Omega$.
\end{thm}
\noindent
\begin{remark}\label{rem:condi} {\rm
 The uniqueness result is conditional in two aspects. First of all, 
 existence of strong solutions clearly cannot be obtained for $d=3$.
 Moreover, for $\delta>0$ (a case which is also not 
 covered by Theorem~\ref{thm: esistenza forte}), we can prove that two
 hypothetical strong solutions emanating from the same initial datum coincide
 only in the case when either of the two satisfies the additional regularity
 \eqref{rego:add}. }
\end{remark}

%%%%%%%%%%%%%%%%%%%%%%%%%%%%%%%%%%%%%%%%%%%%%%%%%%%%%%%%%%%%%%%%%%%%%%%%%%%%%%%%%%%%%%%%%%%%%%%%%%%%%%%%%%%%%%%%%%%%%%%%%%%%%%%%%%%%%%%%%%%%%%%%%%%%%%%

\section{Proofs}
\label{sec:pro}

In this section we will prove Theorems~\ref{thm: esistenza debole}, \ref{thm: esistenza forte}, \ref{thm: unicit} as
well as Proposition~\ref{thm: dissipativita}. We start presenting a number of basic a priori estimates by working
directly, though formally, on system \eqref{eq: N-S}-\eqref{eq: A-C pi_t} without referring explicitly
to any regularization or approximation. Actually, the variational structure underlying system 
\eqref{eq: N-S}-\eqref{eq: A-C pi_t} is rather simple and for this reason it seems to be worth
starting with the computations leading to the energy inequality. This actually constitutes
the main \emph{a priori}\/ information that any (reasonably defined) solution should satisfy. 
Then, a possible approach via a Faedo-Galerkin regularization compatible
with the a priori estimates will be outlined in Section~\ref{sec:Gal app} below.

In what follows we will note as $c$ and $\kc$ some generic positive constants
depending only on the given data of the problem but independent of the final time $T$, with $\kc$ 
being used in estimates from below. The values of $c$ and $\kc$ 
will be allowed to vary on occurrence, and will be assumed to be independent of any hypothetical
regularization or approximation parameter. Specific values of $\kc$ and $c$ will be noted as $\kc_i, c_i$, $i \ge 1$.

%%%%%%%%%%%%%%%%%%%%%%%%%%%%%%%%%%%%%%%%%%%%%%%%%%%%%%%%%%%%%%%%%%%%%%%%%%%%%%%%%%%%%%%%%%%%%%%%%%%%%%%%%%%%%%%%%%%%%%%%%%%%%%%%%%%%%%%%%%%%%%%%%%%%%%%

\subsection{Energy estimate}
\label{subsec:en}

We start testing \eqref{eq: N-S} with $\u$ to obtain
\begin{equation}
\label{eq: 02}
  \frac12 \derivt\normHd{\u}^2 + \normHd{\grad\u}^2 - \ints{(\grad\vp\otimes\grad\vp):\grad\u} = 0,
\end{equation}
where we also used the incompressibility constraint \eqref{eq: N-S div} and the periodic boundary
conditions.
%
%\begin{equation*}
%\ints{\dive(\u\otimes\u)\cdot\u} = \ints{\dive \u \u^2} + b(\u,\u,\u)=0
%\end{equation*}
%thanks to divergence free condition and properties of the trilinear operator $b$ defined in \eqref{def: trilinear b}.  
Next, we test \eqref{eq: A-C pi_t} with $\pi$. Using again periodicity and incompressibility we deduce
\begin{equation}
\label{eq: 01}
  \frac{\ka}{2} \derivt\normH{\pi}^2 + \normH{\pi}^2 + \ints{\grad\pi\cdot\grad\vp} + \ints{f(\vp)\pi} = 0.
\end{equation}
%
%where we used that
%
%\begin{equation*}
%  \delta \ints{(\u\cdot\grad\pi)\cdot\pi}
%   =\frac\delta2\ints{\u \cdot\grad\pi^2} = -\frac\delta2\ints{\pi^2\dive \u}=0.
%\end{equation*}
%
Now we exploit the expression of $\pi$ \eqref{eq: A-C pi} in order to compute the last two 
integrals in~\eqref{eq: 01}. Firstly we have
\begin{equation}
\label{giu:11}
  \ints{\grad\pi\cdot\grad\vp} 
   = \frac12 \ddt \| \nabla \fhi \|_H^2 
   + \sigma \normH{\grad\vp}^2
   + \epsilon \| \Delta \fhi \|_H^2
   + \ints{\nabla \fhi \cdot \nabla ( \u \cdot \nabla \fhi )},
\end{equation}
whereas the second term, thanks also to incompressibility, gives
\begin{equation}
\label{giu:12}
  \ints{f(\vp)\pi} 
   = \ddt \ints{ F(\fhi) }
    + \epsilon \ints{f'(\fhi) | \nabla \fhi |^2 }
    + \sigma \ints{ f(\vp)\vp }.
\end{equation}
Replacing \eqref{giu:11}-\eqref{giu:12} into \eqref{eq: 01}, we deduce
\begin{align}\no
  & \ddt \left( \frac{\ka}{2} \normH{\pi}^2 
     + \frac12 \| \nabla \fhi \|_H^2
     + \ints{ F(\fhi) } \right)
    + \normH{\pi}^2 
    + \sigma\normH{\grad\vp}^2
    + \epsilon \| \Delta \fhi \|_H^2\\ 
 \label{giu:13}
  & \mbox{}~~~~~
   + \sigma \ints{f(\vp)\vp}
   + \ints{\nabla \fhi \cdot \nabla ( \u \cdot \nabla \fhi )}
   + \epsilon \ints{f'(\fhi) | \nabla \fhi |^2 }
    = 0.
\end{align}
We now notice that, by standard integration by parts,
\begin{equation}\label{giu:14}
  \ints{\nabla \fhi \cdot \nabla ( \u \cdot \nabla \fhi )}
   = \ints{(\nabla \fhi \otimes \nabla \fhi) : \nabla \u}
   + \frac12 \ints{\vu \cdot \nabla |\nabla \fhi|^2},
\end{equation}
and the second summand on the right-hand side vanishes by incompressibility. 
Then, summing \eqref{eq: 02} together with \eqref{giu:13} and taking \eqref{giu:14}
into account, we arrive at the \emph{energy estimate}
\begin{align}\no
  & \ddt \left( \frac12 \normHd{\u}^2
     + \frac{\ka}{2} \normH{\pi}^2 
     + \frac12 \| \nabla \fhi \|_H^2
     + \ints{ F(\fhi) } \right)
     + \normHd{\grad\u}^2
     + \normH{\pi}^2    \\ 
 \label{eq:energy}
  & \mbox{}~~~~~
    + \sigma \normH{\grad\vp}^2
    + \epsilon \| \Delta \fhi \|_H^2
    + \sigma \ints{f(\vp)\vp}
    + \epsilon \ints{f'(\fhi) | \nabla \fhi |^2 }
    = 0.
\end{align}
Now, since $F$ is just assumed to be {\sl $\lambda$-convex}, in general
the last term on the \lhs\ may be nonpositive. In order to control it,
we need to exploit the Gronwall inequality. Recalling \eqref{frr} and 
Assumption~\eqref{powerlike} (and using in particular the 
$\lambda$-convexity of $F$), taking a generic $C>0$, we actually deduce
\begin{align}\label{eq: 04}
  & \derivt\bigg[\normHd{\u}^2 + \ka \normH{\pi}^2 + \normH{\grad\vp}^2 + \ints{2F(\vp)} + C\bigg] + 2 \normHd{\grad\u}^2 + 2\normH{\pi}^2  \\
 \nonumber
  & \mbox{}~~~~~ + 2(\sigma-\lambda\epsilon)\normH{\grad\vp}^2+ 2 \epsilon \normH{\lap\vp}^2 
   + 2\sigma\kc_p \bigg( \|\vp\|_{L^{p+2}(\Omega)}^{p+2} + \ints{F(\vp)}\bigg) 
    \le 2\sigma c_p
\end{align}  
and the coefficient $2(\sigma-\lambda\epsilon)$ may be actually nonpositive (and it is
always nonpositive when $\sigma = 0$ and $\lambda>0$).

Let us then note as $\calG$ the sum of the terms in 
square brackets in \eqref{eq: 04}. Now, recalling that $\u$ has zero spatial mean 
by \eqref{hp:init} and \eqref{cons:u}, the Poincar\'e-Wirtinger inequality implies
\begin{equation}\label{eq:powi}
  \normHd{\nabla\u}^2 
   \ge \kc \normHd{\u}^2
\end{equation}  
for some $\kc>0$. Hence, rearranging terms, 
\eqref{eq: 04} gives
\begin{align}\no
  & \ddt \calG 
   + \normHd{\grad\u}^2 
   + \kc \normHd{\u}^2 
   + 2 \normH{\pi}^2 
   + \sigma\kc_p \normH{\grad\vp}^2
   + 2 \sigma\kc_p \ints{F(\vp)} \\
 \label{eq: 04b}
  & \mbox{}~~~~~   
   + 2 \epsilon \normH{\lap\vp}^2
   + 2 \sigma\kc_p \|\vp\|_{L^{p+2}(\Omega)}^{p+2} 
     \le 2\sigma c_p + (\sigma\kc_p - 2\sigma + 2\lambda\epsilon)\normH{\grad\vp}^2.
\end{align}  
Now, using~\eqref{powerlike3} and \eqref{frr2}, it is clear that $C$ can be 
chosen in such a way that
\begin{equation}\label{eq: F bdd}
  \normV{\vp}^2 + k'' \| \vp \|_{L^{p+2}(\Omega)}^{p+2}
   \le \normH{\grad\vp}^2 + 2\ints{F(\vp)} + C 
   \le \normV{\vp}^2 + 2 K'' \| \vp \|_{L^{p+2}(\Omega)}^{p+2} + c.
\end{equation}
%
%for suitable constants $\kappa_p, C_p >0$. 
%
Hence, \eqref{eq: 04b} implies the differential inequality
\begin{equation}\label{eq: 04b2}
  \ddt \calG 
   + \kc \calG 
   + \normHd{\grad\u}^2 
   + 2 \epsilon \normH{\lap\vp}^2
   + 2 \sigma\kc_p \|\vp\|_{L^{p+2}(\Omega)}^{p+2} 
   \le c \big( 1 + \normH{\grad\vp}^2 \big).
\end{equation}  
for some computable constants $\kc,c>0$ depending on the various parameters
of the problem.
Now, when $\sigma=0$ the \rhs\ needs to be estimated by using
Gronwall. Indeed, in that case \eqref{eq: 04b2}, using \eqref{eq: F bdd},
may be rewritten as
\begin{equation}\label{eq: 04b3}
  \ddt \calG 
   + \normHd{\grad\u}^2 
   + 2 \epsilon \normH{\lap\vp}^2
   \le c_1 \calG + c_2,
\end{equation}  
for some constants $c_1,c_2>0$, whence Gronwall's inequality entails
\begin{equation}\label{s:energia1}
 \calG(t) \le \calG(0) e^{c_1 t} 
   + \frac{c_2}{c_1}
\end{equation}
as well as
\begin{equation}\label{s:energia1b}
  \int_0^t \big( \normHd{\grad\u(s)}^2 
   + 2 \epsilon \normH{\lap\vp(s)}^2 \big)\,\dis
  \le c \big( \calG(0) e^{c_1 t} + t \big).
\end{equation}
On the other hand, when $\sigma>0$ it is clear that we can avoid exponential
growth of the energy. Indeed, in that case the \rhs\ of \eqref{eq: 04b2} can be
controlled in the following way: 
\begin{equation}\label{stima: rhs}
  c \normHd{\grad\vp}^2 
    = - c \ints{\vp \lap\vp} 
   \le \epsilon \normH{\lap\vp}^2 + \frac{1}{4\epsilon}\normH{\vp}^2 
    \le \epsilon \normH{\lap\vp}^2 + \sigma\kc_p \norm{\vp}{L^{p+2}}^{p+2} + c(\sigma,\epsilon,p),
\end{equation}
where in the last computation we have used \eqref{frr2}. As a consequence, 
in place of \eqref{eq: 04b3} we get from 
\eqref{eq: 04b2} the better relation
\begin{equation}\label{eq: 04b4}
  \ddt \calG 
   + \kc_1 \calG 
   + \normHd{\grad\u}^2 
   + \epsilon \normH{\lap\vp}^2
   + \sigma\kc_p \|\vp\|_{L^{p+2}(\Omega)}^{p+2} 
   \le c_3,
\end{equation}  
whence we immediately deduce 
\begin{equation}\label{s:energia1c}
  \calG(t) \le \calG(0) e^{-\kc_1 t} 
   + \frac{c_3}{\kc_1}
\end{equation}
and, integrating \eqref{eq: 04b4} over the generic interval $(t,t+1)$, $t\ge 0$,
\begin{equation}\label{dissfiniti}
  \int_{t}^{t+1} 
   \big( \normHd{\grad\u(s)}^2 
      + \epsilon \normH{\lap\vp(s)}^2 
      + \sigma\kc_p \norm{\vp(s)}{L^{p+2}(\Omega)}^{p+2} \big) \, \dis
   \le \calG(0) e^{-\kc_1 t} 
   + c_3 + \frac{c_3}{\kc_1}.
\end{equation}

%%%%%%%%%%%%%%%%%%%%%%%%%%%%%%%%%%%%%%%%%%%%%%%%%%%%%%%%%%%%%%%%%%%%%%%%%%%%%%%%%%%%%%%%%%%%%%%%%%%%%%%%%%%%%%%%%%%%%%%%%%%%%%%%%%%%%%%%%%%%%%%%%%%%%%%

\subsection{Dissipative estimate for $\delta = \ka $}
\label{subsec:diss}

We start with testing \eqref{eq: A-C pi} by $\pi$ to obtain
\begin{equation}\label{eq: d0}
  (\pi,\vp_t) = \normH{\pi}^2 - (\u\cdot\grad\vp,\pi) + \epsilon (\lap\vp,\pi)
   - \sigma (\vp,\pi).
\end{equation}
Next, {\sl assuming $\delta=\kappa$}, we test \eqref{eq: A-C pi_t} by $\vp$. Integrating
by parts in time, we obtain
\begin{equation}\label{eq: d1}
  \ka (\pi,\vp_t) = \ka \ddt (\pi,\vp) + (\pi,\vp) + \ka (\u\cdot\grad\pi,\vp) + \normH{\grad\vp}^2 + (f(\vp),\vp).
\end{equation}
Combining the above relations, we deduce
\begin{align}\no
  & \ka \ddt (\pi,\vp) + (\pi,\vp) 
   + \normH{\grad\vp}^2 + (f(\vp),\vp)\\
 \label{eq: d1b}
   & \mbox{}~~~~~
    = \ka \normH{\pi}^2 - \ka \big[(\u\cdot\grad\vp,\pi) - (\u\cdot\grad\pi,\vp)\big] + \ka \epsilon (\lap\vp,\pi)
    - \kappa \sigma (\vp,\pi).
\end{align}
Now, we observe that, due to $\delta=\kappa$ and incompressibility, 
the transport terms actually vanish:
\begin{equation*}
  (\u\cdot\grad\vp,\pi) + (\u\cdot\grad\pi,\vp) = \ints{\u\cdot\grad(\vp\pi)} 
   = - \ints{ \vp \pi \dive \u} = 0.
\end{equation*} 
Hence relation \eqref{eq: d1b} reduces to
\begin{equation}\label{co:en2}
  \ka \ddt (\pi,\vp) + (\pi,\vp) 
   + \normH{\grad\vp}^2 + (f(\vp),\vp)
  = \ka \normH{\pi}^2 + \ka \epsilon (\lap\vp,\pi)
  - \kappa \sigma (\vp,\pi).
\end{equation}
Let us now multiply the above relation by $\eta > 0$ and sum the result
to~\eqref{eq: 04}. We get
\begin{align}\nonumber
  & \ddt \bigg[ \normHd{\u}^2 + \normH{\pi}^2 + \normH{\grad\vp}^2 + \ints{2F(\vp)} + C + \eta \ka (\pi,\vp) \bigg] \\
 \nonumber
  & \mbox{}~~~~~ + 2 \normHd{\grad\u}^2 + ( 2 - \eta \ka ) \normH{\pi}^2 + 2 \epsilon \normH{\lap\vp}^2 + \eta \normH{\grad\vp}^2
   + \eta (f(\vp),\vp)\\
 \label{eq:04b}  
   & \le 2 \lambda\epsilon\normH{\grad\vp}^2 + \eta \ka \epsilon (\lap\vp,\pi) - \eta(1+\kappa\sigma) ( \pi, \fhi )
\end{align}  
and we need to control the terms on the right-hand side. First of all, we observe that 
\begin{equation}\label{eq:04c}
  \eta\ka\epsilon (\lap\vp,\pi) - \eta(1+\kappa\sigma) ( \pi, \fhi ) 
   \le \frac{\eta \ka^2 \epsilon^2}2 \normH{\lap\vp}^2
    + \frac{\eta}2 \norm{\fhi}{H}^2  
    + \frac{\eta}{2} \big(1+(1+\ka\sigma)^2\big) \norm{\pi}{H}^2.
\end{equation}
On the other hand,
\begin{equation}\label{eq:04d}
  2 \lambda\epsilon\normH{\grad\vp}^2
   = 2 \lambda\epsilon \io ( - \Delta\fhi ) \fhi \, \dix
   \le \frac\epsilon2 \normH{\lap\vp}^2
    + 2\lambda^2 \epsilon \normH{\vp}^2.
\end{equation}
Replacing \eqref{eq:04c}-\eqref{eq:04d} into \eqref{eq:04b} we then deduce
\begin{align}\nonumber
  & \ddt \bigg[ \normHd{\u}^2 + \normH{\pi}^2 + \normH{\grad\vp}^2 + \ints{2F(\vp)} + C + \eta \ka (\pi,\vp) \bigg] 
     + 2 \normHd{\grad\u}^2 \\
 \label{eq:04e}  
  & \mbox{}~~~~~ + \Big( 2 - \eta\kappa - \frac\eta2\big(1+(1+\ka\sigma)^2\big)  \Big) \normH{\pi}^2 
   + \Big( \frac{3\epsilon}2 - \frac{\eta \ka^2 \epsilon^2}2 \Big) \normH{\lap\vp}^2 
   + \eta \normH{\grad\vp}^2
   + \eta (f(\vp),\vp)\\ \nonumber 
  & \mbox{}~~~~~ \le \Big( \frac{\eta}2 + 2 \lambda^2 \epsilon \Big) \normH{\vp}^2.
\end{align}  
Now, we specify the choice of the ``small'' positive constant $\eta$. 
First, we need $\kappa\eta \le 2/3$, in such a way that
\begin{align}\label{eq:04d2}
  \kappa\eta (\pi,\vp)
  & \ge - \frac{1}{3} \| \pi \|^2_H - \frac{1}{3} \| \vp \|^2_H
  \ge - \frac{1}{3} \| \pi \|^2_H - \io F(\fhi)\,\dix - c,\\
 \label{eq:04d3}
  \kappa\eta (\pi,\vp)
  & \le \frac{1}{3} \| \pi \|^2_H + \frac{1}{3} \| \vp \|^2_H,
\end{align}
for some $c > 0$, where also \eqref{frr2} has been used. Second, we need $\eta$ 
so small that the coefficients of the first two terms on the second row of 
\eqref{eq:04e} are strictly positive. We also note that, for such $\eta$,
using $\calD$ for the sum of the terms in 
square brackets in \eqref{eq:04e}, there exists 
(a new value of) $C>0$ such that
\begin{align}
 \label{eq:calD}     
  & \normH{\grad \vp}^2
  + \frac{2}{3} \normH{\pi}^2
  + \normHd{\u}^2
  + \ints{F(\vp)}
   \le \calD \\
  &\mbox{}~~~~~  \le \frac{4}{3} \normV{\vp}^2
   + \frac43 \normH{\pi}^2 
   + \normHd{\u}^2 
   + \ints{2F(\vp)} 
   + C.
\end{align}
Next, thanks to \eqref{frr} and \eqref{frr2}, the 
right-hand side of \eqref{eq:04e} can be estimated as follows
\begin{equation}\label{eq:04e2a}
  \Big( \frac{\eta}2 + 2 \lambda^2 \epsilon \Big) \normH{\vp}^2
   \le \frac\eta2 (f(\vp),\vp) + c(\eta,\epsilon,\lambda),
\end{equation}  
and the nonconstant term can be absorbed by the corresponding one
on the left-hand side. Moreover, the term depending on $f$ on the \lhs\
can be estimated from below by means of \eqref{frr}.
Hence, using once more the Poincar\'{e}-Wirtinger inequality for $\u$, 
we see that \eqref{eq:04e} implies, for some $\kc_2,\kc_3>0$, the 
differential inequality
\begin{equation}\label{eq:diss}
  \ddt \calD 
   + \kc_2 \calD 
   + \kc_3 \normHd{\nabla\u}^2 
   + \kc_3 \epsilon \normH{\lap\vp}^2
    \le c_4 ,
\end{equation}  
with the constants $\kc_2, \kc_3,c_4>0$ depending on $\epsilon,\lambda$ and also through
the choice of $\eta$. As a consequence, we deduce
\begin{equation}\label{diss:est}
  \calD(t) 
   + \int_t^{t+1} \kc_3 \big( \normHd{\nabla\u(s)}^2 + \epsilon \normH{\lap\vp(s)}^2 \big) \,\dis
    \le \calD(0) e^{-\kc_2 t} + \frac{c_4}{\kc_2} + c_4 ,
    \quad \perogni t \ge 0.   
\end{equation}

%%%%%%%%%%%%%%%%%%%%%%%%%%%%%%%%%%%%%%%%%%%%%%%%%%%%%%%%%%%%%%%%%%%%%%%%%%%%%%%%%%%%%%%%%%%%%%%%%%%%%%%%%%%%%%%%%%%%%%%%%%%%%%%%%%%%%%%%%%%%%%%%%%%%%%%

\subsection{Proof of Theorem \ref{thm: esistenza debole}: existence of weak solutions}
\label{subsec:wss}

In the sequel we consider a sequence $\{(\u_n,\fhi_n,\pi_n)\}$, $n\in\NN$,
of triplets satisfying system \eqref{eq: N-S}-\eqref{eq: A-C pi_t} in 
a suitable sense and complying, uniformly with respect to $n$,
with the energy principle. More precisely, we assume that
relation~\eqref{eq: 04b2} is satisfied with constants $k, c$ independent of
$n$. Hence, either \eqref{s:energia1}-\eqref{s:energia1b} 
or \eqref{s:energia1c}-\eqref{dissfiniti} also hold 
with constants independent of $n$. As a consequence, we can prove that 
any limit point (in a suitable sense) of a generic (nonrelabelled) subsequence of 
$\{(\u_n,\fhi_n,\pi_n)\}$ is a weak solution in the sense of Theorem~\ref{thm: esistenza debole}.
In Section~\ref{sec:Gal app} we will see how the present argument can be
adapted to a Galerkin approximation scheme.

That said, we first observe that, from \eqref{s:energia1}-\eqref{s:energia1b} 
(or from \eqref{s:energia1c}-\eqref{dissfiniti}), using the definition of $\calG$
(cf.~\eqref{eq: 04}) and the structure property \eqref{eq: F bdd}, 
we may deduce the following relations:
\begin{align} \label{conv:11}     
  & \u_n \to \u \quext{weakly star in }\/ L^\infty(0,T;\H) \cap L^2(0,T;\V),\\
 \label{conv:12}     
  & \fhi_n \to \fhi \quext{weakly star in }\/ L^\infty(0,T;V) \cap L^\infty(0,T;L^{p+2}(\Omega)) 
      \cap L^2(0,T;H^2(\Omega)),\\
 \label{conv:13}     
  & \pi_n \to \pi \quext{weakly star in }\/ L^\infty(0,T;H).
\end{align}
Here and below all the convergence properties will be intended to hold
up to the extraction of (nonrelabelled) subsequences of $n\nearrow\infty$.

We now deduce some further consequence of the energy estimate that will be
needed in order to take the limit of the product terms in the system. 
We start dealing with the velocity and observe that, for any 
$\v\in \V$,
\begin{equation} \label{co:001}
  \ints{ (\u_n \cdot \nabla \u_n) \cdot \v }
   = - \ints{ (\u_n \otimes \u_n) : \nabla \v }
   \le \| \u_n \|_{L^4(\Omega)}^2 \| \nabla \v \|_H.
\end{equation}
Now, for $d=2$, thanks to Ladyzhenskaya's inequality, \eqref{conv:11} implies 
\begin{equation} \label{co:002}
  \| \u_n \|_{L^4(0,T;L^4(\Omega))} \le c,
\end{equation}
whence from \eqref{co:001} we obtain
\begin{equation} \label{co:003}
  \| \u_n \cdot \nabla \u_n \|_{L^2(0,T;\V')} \le c.
\end{equation}
On the other hand, for $d=3$, Sobolev's embeddings and interpolation only give
\begin{equation} \label{co:002b}
  \| \u_n \|_{L^{8/3}(0,T;L^4(\Omega))} \le c,
\end{equation}
whence from \eqref{co:001} we deduce
\begin{equation} \label{co:003b}
  \| \u_n \cdot \nabla \u_n \|_{L^{4/3}(0,T;\V')} \le c.
\end{equation}
Analogously, from \eqref{conv:12} we have
\begin{align} \label{co:002x}
  & \| \nabla \fhi_n \|_{L^4(0,T;L^4(\Omega))} \le c 
   \quext{for }\,d=2,\\
 \label{co:003x}
  & \| \nabla \fhi_n \|_{L^{8/3}(0,T;L^4(\Omega))}
   % + \| \nabla \fhi_n \|_{L^{4}(0,T;L^3(\Omega))} 
    \le c 
   \quext{for }\,d=3,
\end{align}
whence 
\begin{align} \label{co:002y}
  & \| \nabla \fhi_n \otimes \nabla \fhi_n \|_{L^2(0,T;H)} \le c 
   \quext{for }\,d=2,\\
 \label{co:003y}
  & \| \nabla \fhi_n \otimes \nabla \fhi_n \|_{L^{4/3}(0,T;H)} \le c 
   \quext{for }\,d=3.
\end{align}
Now, testing \eqref{eq: N-S} (written for $\u_n,\fhi_n$) by $\v\in\V$ and rearranging,
we have
\begin{equation}\label{eq: 001}
  \duav{\u_{n,t},\v} = - \duav{\nabla\u_n,\nabla \v} 
   + \duav{ \grad\vp_n\otimes\grad\vp_n, \nabla \v } 
   + \duav{ \u_n \otimes \u_n, \nabla \v },
\end{equation}
whence, using \eqref{co:003}, \eqref{co:003b}
and \eqref{co:002y}-\eqref{co:003y}, we easily infer
\begin{equation} \label{conv:15}     
   \u_{n,t} \to \u_t \quext{weakly in }\/ L^q(0,T;\V'),
\end{equation}
where $q=2$ if $d=2$ and $q=4/3$ if $d=3$. Combining the above with 
\eqref{conv:11} and applying the Aubin-Lions lemma, we then obtain
\begin{equation} \label{conv:16}     
   \u_{n} \to \u \quext{strongly in }\/ C^0([0,T];\V') \cap L^2(0,T;\H).
\end{equation}

\smallskip

We now move to considering the behavior of $\fhi_n$. To this aim, we first
notice that, using \eqref{co:002} and \eqref{co:002x} in 2D,
and using the first \eqref{conv:11}, the last \eqref{conv:12} and Sobolev's
embeddings in 3D, there follows
\begin{equation} \label{conv:17}
   \| \u_{n}\cdot \nabla \fhi_n \|_{L^2(0,T;L^s(\Omega))} \le c,
\end{equation}
where $s=3/2$ if $d=3$ and $s =2$ if $d=2$. Hence, comparing 
terms in \eqref{eq: A-C pi} and using again \eqref{conv:11}-\eqref{conv:13}
it is not difficult to deduce
\begin{equation} \label{conv:21}     
   \vp_{n,t} \to \vp_t \quext{weakly in }\/ L^2(0,T;L^s(\Omega)),
\end{equation}
with $s$ as above. Applying once more the Aubin-Lions lemma, we 
infer
\begin{equation} \label{conv:22}     
   \vp_{n} \to \vp \quext{strongly in }\/ C^0([0,T];H) \cap L^2(0,T;V).
\end{equation}
As a consequence, we obtain (at least) that
\begin{equation} \label{conv:23}     
   \nabla\vp_{n} \otimes \nabla\vp_n \to \nabla\vp \otimes \nabla\vp
    \quext{strongly in }\/ L^1(0,T;L^1(\Omega)).
\end{equation}
This allows us to take the limit $n\nearrow \infty$ in relation \eqref{eq: 001}
to get back \eqref{N-S:weak}. Moreover, it is clear that, from \eqref{conv:16}
and \eqref{conv:22}, there also follows (at least)
\begin{equation} \label{conv:24}     
   \u_n\cdot \nabla \vp_{n} \to \u \cdot \nabla \vp 
    \quext{strongly in }\/ L^1(0,T;L^1(\Omega)).
\end{equation}
Actually, relations \eqref{conv:23} and \eqref{conv:24} could be improved,
anyway they suffice for the sequel. Indeed, we can now 
write \eqref{eq: A-C pi} at the $n$-level, take 
the limit $n\nearrow \infty$, and get back \eqref{AC:pi:weak}.

Finally, we need to pass to the limit in \eqref{eq: A-C pi_t}. Then, let us 
notice that \eqref{conv:12} and \eqref{powerlike2} entail
\begin{equation} \label{conv:25}     
   \| f(\vp_n) \|_{L^\infty(0,T;L^{\frac{p+2}{p+1}}(\Omega))} \le c. 
\end{equation}
Combining this fact with the pointwise convergence resulting from
\eqref{conv:22} we then deduce
\begin{equation} \label{conv:26}     
  f(\vp_n) \to f(\vp) \quext{weakly star in }\, L^\infty(0,T;L^{\frac{p+2}{p+1}}(\Omega))
   \quext{and strongly in }\, L^1(0,T;L^{1}(\Omega)).
\end{equation}
Finally, if $\delta>0$, we need to take care of the convection term
in~\eqref{eq: A-C pi_t}. Combining \eqref{conv:11} with \eqref{conv:13} and 
using interpolation with the 2D embedding $H^1(\Omega) \immcont L^q(\Omega)$ for 
$1 \le q < \infty$, we actually deduce that
\begin{equation} \label{conv:27}     
  \| \u_n \pi_n \|_{L^2(0,T;L^r(\Omega))} \le c_r,
\end{equation}
where $r = 3/2$ for $d=3$ and $r\in[1,2)$ for $d=2$. Then, using
the second of \eqref{conv:16}, we infer
\begin{equation} \label{conv:28}     
  \u_n \pi_n \to \u \pi \quext{weakly in }\,L^2(0,T;L^r(\Omega)).
\end{equation}
Hence, for any $v\in W^{1,3}(\Omega)$, using incompressibility and 
periodic boundary conditions, we deduce 
\begin{equation} \label{conv:29}     
  \ints{ (\u_n \cdot \nabla \pi_n) v } 
    = - \ints{ \pi_n (\u_n \cdot \nabla v) } 
    \to - \ints{ \pi (\u \cdot \nabla v) }
   \quext{weakly in }\, L^2(0,T)
\end{equation}
and in particular the distributional divergence
$\dive ( \pi \u )$ lies in the dual space $L^2(0,T;(W^{1,3})'(\Omega))$.

Using \eqref{conv:11}-\eqref{conv:13} together with
\eqref{conv:25} and \eqref{conv:27}, comparing terms in the $n$-version
of \eqref{eq: A-C pi_t}, and using continuity of the embedding
$L^\infty(0,T;L^{\frac{p+2}{p+1}}(\Omega)) \subset L^2(0,T;(W^{1,3})'(\Omega))$,
we then obtain
\begin{equation} \label{conv:30}
  \pi_{n,t} \to \pi_t \quext{weakly in }\, 
    L^2(0,T;(W^{1,3})'(\Omega)).   
\end{equation}
This fact, by the Aubin-Lions lemma and \eqref{conv:13}, also yields
the strong convergence
\begin{equation} \label{conv:31}
  \pi_{n} \to \pi \quext{strongly in }\, 
    C^0([0,T];V').
\end{equation}
The above relations permit us to
pass to the limit in \eqref{eq: A-C pi_t}, which may be reinterpreted as an equality
in $L^2(0,T;(W^{1,3})'(\Omega))$ (or also as an equality in $(W^{1,3})'(\Omega)$,
almost everywhere in time). Namely, in the limit we get back exactly \eqref{AC:pi_t:weak}.
Note also that the regularity conditions \eqref{rego:u:w}-\eqref{rego: pi t}
in the statement of Theorem~\ref{thm: esistenza debole}
are a direct consequence of \eqref{conv:11}-\eqref{conv:13}, \eqref{conv:15}
\eqref{conv:21} and \eqref{conv:30}.

Finally, we observe that the strong convergence relations \eqref{conv:16}, \eqref{conv:22},
\eqref{conv:31} imply
\begin{equation} \label{conv:32}
  (\u_n,\vp_n,\pi_n)|_{t=0} 
   \to (\u,\vp,\pi)|_{t=0} \quext{strongly in }\, 
   \V'\times H \times V'.
\end{equation}
Hence, assuming that the triplet $(\u_n,\vp_n,\pi_n)$ satisfies an initial
condition of the form
\begin{equation} \label{conv:33}
  (\u_n,\vp_n,\pi_n)|_{t=0} 
   = (\u_{0,n},\vp_{0,n},\pi_{0,n}),
\end{equation}
where $(\u_{0,n},\vp_{0,n},\pi_{0,n})$ tends to $(\u_{0},\vp_{0},\pi_{0})$ in 
a suitable way, letting $n\nearrow \infty$, we obtain \eqref{init} in the limit,
which concludes the proof of Theorem~\ref{thm: esistenza debole}.
\begin{remark}\label{rem:eps}{\rm
 It is worth discussing a bit more the occurrence of the viscous regularization
 term in \eqref{eq: A-C pi} at the light of the a priori estimates. Actually,
 if that term is omitted (i.e., if $\epsilon=0$), the latter convergence in
 \eqref{conv:12} would be lacking and we would not have any $L^p$-information on second space
 derivatives of $\fhi$. Hence, the $L^2(0,T;V)$-convergence in \eqref{conv:22}
 would also be missing and we could not take the limit of the \rhs\ of \eqref{eq: N-S} 
 as specified in \eqref{conv:23}.
 }
\end{remark}
%

%%%%%%%%%%%%%%%%%%%%%%%%%%%%%%%%%%%%%%%%%%%%%%%%%%%%%%%%%%%%%%%%%%%%%%%%%%%%%%%%%%%%%%%%%%%%%%%%%%%%%%%%%%%%%%%%%%%%%%%%%%%%%%

\subsection{Proof of Proposition~\ref{thm: dissipativita}: dissipativity}
\label{ssec:diss}

In the case $\sigma>0$ we can take advantage of estimate \eqref{s:energia1c},
whereas for $\delta=\ka$ we have relation \eqref{diss:est}. Of course,
these bounds hold a priori for the approximating solutions $(\u_n,\vp_n,\pi_n)$.
On the other hand, they pass to the limit $n\nearrow\infty$ because the quantities
on the \rhs s are independent of $n$ and we can use semicontinuity of norms 
with respect to weak or weak star convergence when we take the limit. Hence we obtain
that \eqref{s:energia1c}, or \eqref{diss:est}, also holds for 
limit solution(s) $(\u,\vp,\pi)$. Noting that the functionals $\calG$ (cf.~\eqref{eq: 04})
and $\calD$ (cf.~\eqref{eq:04e}) control both from above and from below
the norms specified in the statement, it is then apparent that
\eqref{s:energia1c} (or \eqref{diss:est}) implies the desired bound 
\eqref{st:diss} for every $t\ge T_0$, with $T_0$ depending only on the 
magnitude of the initial energy, i.e., on the norms in~\eqref{magn:init}.
\begin{remark}\label{rem:attr}{\rm
 One may wonder if the dissipative estimate \eqref{st:diss}
 (possibly combined with the regularity estimates obtained in the proof of 
 existence of strong solutions) could be used to prove existence 
 of the global attractor at least in the 2D case. We do not address 
 this interesting issue here, but we limit ourselves to observe that 
 the problem seems nontrivial and its resolution may require the 
 use of some careful decomposition method. Indeed, 
 while the variables $\u$ and $\fhi$ enjoy some 
 regularization property (if they are, respectively, 
 in $\H$ and in $V$ at the initial time, then they 
 are in $\V$ and in $H^2(\Omega)$ 
 for some small $t>0$), i.e., they have a {\sl parabolic}\/ behavior, 
 on the other hand the variable $\pi$ has a {\sl hyperbolic}\/ behavior,
 i.e.~it does not seem to regularize in time (compare \eqref{rego:pi:w} for $\pi$
 with \eqref{rego:u:w} and \eqref{rego:fhi:w} for $\u$ and $\fhi$). 
 }
\end{remark}

%%%%%%%%%%%%%%%%%%%%%%%%%%%%%%%%%%%%%%%%%%%%%%%%%%%%%%%%%%%%%%%%%%%%%%%%%%%%%%%%%%%%%%%%%%%%%%%%%%%%%%%%%%%%%%%%%%%%%%%%%%%%%%

\subsection{Proof of Theorem~\ref{thm: esistenza forte}: existence of strong solutions}
\label{ssec:sol for}

The proof is based on some additional regularity estimates. As before, we work directly on system
\eqref{eq: N-S}-\eqref{eq: A-C pi_t}. We postpone to the next section
the justification of this argument at the light of the approximation scheme.

First of all, we multiply \eqref{eq: N-S} with $-\lap\u$ and integrate over $\Omega$ to obtain
\begin{align} \label{eq: 1}
  0 & = \int_{\Omega} \big(-\u_t \lap\u + |\lap\u|^2 - \dive(\grad\varphi\otimes\grad\varphi)\lap\u 
    - \grad p \lap \u - \dive(\u\otimes\u)\lap \u \big) \dx\\
 \nonumber
    & = \frac{1}{2} \ddt \norm{\grad\u}{\H}^2 + \norm{\lap\u}{\H}^2 
     + \int_{\Omega}\big(- (D^2\varphi \grad\varphi)\cdot \lap\u - \lap\varphi(\grad\varphi\cdot \lap\u) \big) \dx.
\end{align} 
Indeed, the convection term on the first row vanishes due to $d=2$ (see, e.g.,
\cite[p.~242]{Robinson}).

Now using Ladyzhenskaya's and Young's inequalities with 2D Sobolev embeddings, 
we estimate the last two terms as follows:
\begin{equation}\label{eq: 1'}
 \begin{aligned}
 & \bigg| \int_{\Omega}\big( - (D^2\varphi \grad\varphi)\cdot \lap\u - \lap\varphi(\grad\varphi\cdot \lap\u) \big) \dx \bigg| \\
  & \mbox{}~~~~~ \le \norm{\grad\varphi}{L^4(\Omega)}\norm{D^2\vp}{L^4(\Omega)}\norm{\lap\u}{\H} \\
  & \mbox{}~~~~~ \le \norm{\grad\varphi}{H}^{1/2} \norm{\grad\varphi}{V}^{1/2} 
       \norm{D^2\varphi}{H}^{1/2} \norm{D^2\varphi}{V}^{1/2} \norm{\lap\u}{\H}\\
  & \mbox{}~~~~~ \le \alpha\normHd{\lap\u}^2 
   + c_\alpha \norm{\grad\varphi}{V} \norm{D^2\varphi}{H} \norm{D^3\varphi}{H}\\
  & \mbox{}~~~~~ \le \alpha\normHd{\lap\u}^2 + \alpha\norm{D^3\varphi}{H}^2
   + c_\alpha \norm{D^2\varphi}{H}^4\\ 
  & \mbox{}~~~~~ \le \alpha\normHd{\lap\u}^2 + \alpha\norm{\nabla\lap \varphi}{H}^2
   + c_\alpha \norm{\Delta\varphi}{H}^4, 
 \end{aligned}
\end{equation}
where $\alpha>0$ denotes small constants to be chosen later and
the constants $c_\alpha>0$ are correspondingly
large. In the above computation we used the uniform in time $V$-bound for $\fhi$ 
resulting from the energy bound (cf.~\eqref{conv:12}) and elliptic regularity.

\smallskip

Next, we deal with the Allen-Cahn system. Actually, testing \eqref{eq: A-C pi} with $\lap^2\vp$, 
\eqref{eq: A-C pi_t} with $-\lap\pi$, summing the results, and performing standard manipulations
(note in particular that a couple of terms cancels out), we obtain
\begin{equation}\label{eq: 2}
\begin{aligned}
  & \frac12\derivt\bigg[ \ka \norm{\grad\pi}{H}^2 +
   \norm{\lap\varphi}{H}^2\bigg] + \norm{\grad\pi}{H}^2
  + \sigma \normH{\lap\vp}^2
  + \epsilon\norm{\grad\lap\varphi}{H}^2\\
  & \mbox{}~~~~~ 
   = - \ints{f'(\varphi)\grad\varphi\grad\pi} 
    + \ints{(\nabla \u \nabla \fhi) \cdot \grad \lap \fhi}  
    + \ints{(D^2 \fhi \u) \cdot \grad \lap \fhi}  
   =: \sum_{j=1}^3 I_j
\end{aligned}
\end{equation}
and we need to control the terms $I_j$ on the \rhs. We actually have
\begin{align}\label{eq: 2'}
  | I_1 | & \le c \ints{(1+|\varphi|^p) |\grad\varphi||\grad\pi|} \\
 \nonumber
    & \le c \norm{|\varphi|^p}{L^4(\Omega)}\norm{\grad\varphi}{L^4(\Omega)}\norm{\grad\pi}{H} + c \normH{\grad\vp} \normH{\grad\pi}\\
 \nonumber
    & \le c \big( 1 + \norm{\varphi}{L^{4p}(\Omega)}^p \| \grad\fhi \|_{H}^{1/2} \| \grad\fhi \|_{V}^{1/2} \big) \normH{\grad\pi}\\
 \nonumber
    & \le c \big( 1 + \| \fhi \|_{H^2(\Omega)}^{1/2} \big) \normH{\grad\pi}\\
 \nonumber
    & \le \alpha \normH{\grad\pi}^2 + c_{\alpha} \big( 1 + \| \lap\fhi \|_{H} \big),
 \end{align}
 \begin{align}
 \label{eq: 2''}
 | I_2 | & \le \norm{\grad\u}{\H}^{1/2} \norm{\grad\u}{\V}^{1/2} \norm{\grad\varphi}{H}^{1/2} \norm{\grad\varphi}{V}^{1/2} \norm{\grad\lap\varphi}{H}\\
  \nonumber 
   & \le \alpha \norm{\grad\u}{\V}^2  + \alpha \norm{\grad\lap\varphi}{H}^2 +c_{\alpha} \norm{\grad\varphi}{H}^2 \norm{\grad\varphi}{V}^2 \norm{\grad\u}{\H}^2\\
  \nonumber
   & \le \alpha \norm{\Delta\u}{\H}^2  + \alpha \norm{\grad\lap\varphi}{H}^2 + c_{\alpha} \norm{\lap\varphi}{H}^2 \norm{\grad\u}{\H}^2,
 \end{align}
 \begin{align}
 \label{eq: 2'''}
  | I_3 | & \le \norm{\u}{L^4(\Omega)}\norm{D^2\varphi}{L^4(\Omega)}\norm{\grad\lap\varphi}{H}\\
  \nonumber 
   & \le c \norm{\u}{\H}^{1/2} \norm{\u}{\V}^{1/2} \norm{D^2\varphi}{H}^{1/2}\norm{D^3\varphi}{H}^{1/2} \norm{\grad\lap\varphi}{H}\\
  \nonumber 
   & \le c \norm{\u}{\V}^{1/2} \norm{\lap\varphi}{H}^{1/2} \norm{\grad\lap\varphi}{H}^{3/2} \\
 \nonumber
   & \le \alpha \norm{\grad\lap\varphi}{H}^2 + c_\alpha \norm{\nabla\u}{\H}^2 \norm{\lap\varphi}{H}^2,
\end{align}
where $\alpha$ and $c_\alpha$ are as above (of course their value can vary on occurrence) and we have 
repeatedly used assumption \eqref{powerlike},
2D-Sobolev embeddings, Young's and Ladyzhenskaya's inequalities, and the information 
already obtained with the energy bound (cf.~\eqref{conv:11}-\eqref{conv:13}).

\smallskip

Now we sum \eqref{eq: 1} to \eqref{eq: 2}. Using \eqref{eq: 1'}, \eqref{eq: 2'}, \eqref{eq: 2''} and \eqref{eq: 2'''},
we then deduce
\begin{align}\label{eq: stima}
  & \frac12 \derivt \Big[\normHd{\grad\u}^2 + \ka \normH{\grad\pi}^2 + \normH{\lap\vp}^2\Big] 
   + \normHd{\lap\u}^2 + \normH{\grad\pi}^2 + \sigma \normH{\lap\vp}^2 + \epsilon \normH{\grad\lap\vp}^2\\
 \nonumber
  & \mbox{}~~~~~ \le \alpha\normHd{\lap\u}^2 + \alpha\norm{\nabla\lap \varphi}{H}^2
   + c_\alpha + c_\alpha \norm{\Delta\varphi}{H}^4
   + \alpha \normH{\grad\pi}^2
   + c_{\alpha} \norm{\lap\varphi}{H}^2 \norm{\grad\u}{\H}^2.
\end{align}  
Now, we can take the various $\alpha$'s small enough (in a way that also depends on $\epsilon$)
so that the above reduces to 
\begin{align}\label{eq: stima:2}
  & \derivt \Big[\normHd{\grad\u}^2 + \ka \normH{\grad\pi}^2 + \normH{\lap\vp}^2\Big] 
   + \normHd{\lap\u}^2 + \normH{\grad\pi}^2 + \sigma \normH{\lap\vp}^2 + \epsilon \normH{\grad\lap\vp}^2\\
 \nonumber
  & \mbox{}~~~~~ \le c \big( 1 + \normH{\lap\fhi}^2 \big)
  \big( 1 + \norm{\lap\varphi}{H}^2 + \norm{\grad\u}{\H}^2 \big).
\end{align}  
This inequality has the structure
\begin{equation}\label{eq: stima3}
  \ddt \calE_1(t) + \calD_1(t)
   \le m(t) \big( 1 + \calE_1(t) \big)
\end{equation}
where $\calE_1$ denotes the quantity in square brackets on the \lhs, $\calD_1$ is 
the sum of the remaining terms, and the function 
\begin{equation}\label{eq: stima3b}
  t \mapsto m(t):= c \big( 1 + \norm{\lap\varphi}{H}^2 \big)
\end{equation}
belongs to $L^1(0,T)$ thanks to the $L^2(0,T;H^2(\Omega))$-bound for $\fhi$ 
following from the energy estimate (cf.~\eqref{conv:12}). 

Applying Gronwall's Lemma to the function $t\mapsto \calE_1(t)$ we then conclude that
\begin{equation}\label{reg: a1}
 \begin{aligned}
  & \u \in L^\infty(0,T;\V),\\
  & \pi \in L^\infty(0,T;V),\\
  & \varphi \in L^\infty(0,T;H^2(\Omega)).
 \end{aligned}
\end{equation}
Note that the improved regularity assumptions on the initial data
\eqref{hp:init2} have also been used here.

Then, integrating once more \eqref{eq: stima3} in time we also deduce
\begin{equation}\label{reg: a2}
\begin{aligned}
   &\u \in L^2(0,T;H^2(\Omega)),\\
   &\grad\lap\varphi \in L^2(0,T;H),
\end{aligned}
\end{equation}
whence, using once more elliptic regularity we infer
\begin{equation}\label{reg: a3}
  \vp \in L^2(0,T;H^3(\Omega)).
\end{equation}
Finally, comparing terms in the equations of the system and using the above bounds,
standard manipulations permit us to deduce also
\begin{equation}\label{reg: a4}
\begin{aligned}
   &\u_t \in L^2(0,T;\H),\\
   &\fhi_t \in L^2(0,T;V),\\
   &\pi_t \in L^2(0,T;H),
\end{aligned}
\end{equation}
which is the last bound we need.

\smallskip

Now, in order to complete the proof, we interpret the information obtained
above in the framework of an approximation scheme. In this respect,
\eqref{reg: a1}, \eqref{reg: a2}, \eqref{reg: a3}, and \eqref{reg: a4} can be seen
as a priori estimates uniform with respect to the approximation
parameter $n$. Consequently, the convergence relations \eqref{conv:11}-\eqref{conv:13}
may also be improved in the corresponding way. Finally, we may observe that,
in view of the enhanced regularity, all terms in equations \eqref{eq: N-S} and 
\eqref{eq: A-C pi_t} lie in some $L^p$-space and may be consequently
interpreted in the pointwise sense, as noted in the statement of the theorem.
This completes the proof of Theorem~\ref{thm: esistenza forte}.

%%%%%%%%%%%%%%%%%%%%%%%%%%%%%%%%%%%%%%%%%%%%%%%%%%%%%%%%%%%%%%%%%%%%%%%%%%%%%%%%%%%%%%%%%%%%%%%%%%%%%%%%%%%%%%%%%%%%%%%%%%%%%%

%%%%%%%%%%%%%%%%%%%%%%%%%%%%%%%%%%%%%%%%%%%%%%%%%%%%%%%%%%%%%%%%%%%%%%%%%%%%%%%%%%%%%%%%%%%%%%%%%%%%%%%%%%%%%%%%%%%%%%%%%%%%%%

\subsection{Proof of Theorem~\ref{thm: unicit}: uniqueness}
\label{ssec:unicit}

Let us given a pair of strong solutions $(\u_1, \varphi_1, \pi_1)$ and $(\u_2, \varphi_2, \pi_2)$
of \eqref{eq: N-S}-\eqref{eq: A-C pi_t} emanating from the same initial datum $(\u_0, \varphi_0, \pi_0)$
satisfying \eqref{hp:init2}. We will show that the two solutions do coincide. We perform the proof
in the more difficult case $d=3$ (as said, in this situation the result is conditional because
strong solutions are not known to exist); of course the argument extends to $d=2$ where we actually have
better embeddings. That said, we put $\u \ede \u_1 - \u_2$, $\varphi \ede \varphi_1 - \varphi_2$,
$\pi \ede \pi_1 - \pi_2$ and $p := p_1 - p_2$. Then, writing \eqref{eq: N-S} for the two solutions,
and taking the difference, we get
\begin{equation}\label{eq: N-S diff}
  \u_t 
   - \lap\u 
   + \dive(\u \otimes\u_1) + \dive(\u_2 \otimes \u ) 
   + \grad p
   + \dive(\grad\varphi \otimes\grad\varphi_1) + \dive(\grad\varphi_2\otimes\grad\varphi) 
  = 0.
\end{equation}
Proceeding in the same way for \eqref{eq: A-C pi} and \eqref{eq: A-C pi_t} we obtain
\begin{align} \label{eq: A-C pi diff}
  & \pi = \varphi_t 
   + \u_1 \cdot \grad\varphi + \u \cdot \grad\varphi_2 
   - \epsilon\lap \varphi 
   + \sigma \vp,\\
 \label{eq: A-C pi_t diff}
  & \ka \pi_t + \pi 
   + \delta \u_1 \cdot \nabla \pi + \delta \u \cdot \nabla \pi_2 
   - \lap\varphi 
   + f(\varphi_1)- f(\varphi_2) = 0.
\end{align}
We now work on system \eqref{eq: N-S diff}-\eqref{eq: A-C pi_t diff} with the aim
of getting a contractive estimate. 
We start testing \eqref{eq: N-S diff} with $\u$ to get
\begin{align}\label{co:51}
  & \frac12 \derivt\norm{\u}{\H}^2 + \norm{\grad\u}{\H}^2 
   = \ints{  \big( (\grad\varphi \otimes\grad\varphi_1) + (\grad\varphi_2\otimes\grad\varphi) \big): \nabla \u }\\
 \nonumber
   & \mbox{}~~~~~ + \ints{  \big( (\u \otimes\u_1) + (\u_2 \otimes \u )  \big): \nabla \u }
    =: K_1 + K_2.
\end{align}
Let us now provide a bound of the terms $K_i$ on the \rhs. First, we have
\begin{align}\label{st:K1}
  | K_1 | & \le \norm{\grad\u}{\H} \norm{\grad\fhi}{L^3(\Omega)} \big( \norm{\grad\fhi_1}{L^6(\Omega)} + \norm{\grad\fhi_2}{L^6(\Omega)} \big) \\
 \nonumber
   & \le c \norm{\grad\u}{\H} \norm{\grad\fhi}{H}^{1/2} \norm{\lap\fhi}{H}^{1/2} 
   \le \frac18 \norm{\grad\u}{\H}^2 + \frac{\epsilon}6 \norm{\lap\fhi}{H}^{2} + c_\epsilon \norm{\grad\fhi}{H}^2,
\end{align}
where we have used Sobolev's embeddings, Young's inequality, regularity \eqref{rego:fhi:s} both for 
$\fhi_1$ and for $\fhi_2$, and elliptic regularity results (in order to control the 
$V$-norm of $\grad\fhi$ with the $H$-norm of $\lap\fhi$). Secondly, we have
\begin{align}\label{st:K2}
  |K_2| & \le \norm{\grad\u}{\H} \norm{\u}{L^3(\Omega)} \big( \norm{\u_1}{L^6(\Omega)} + \norm{\u_2}{L^6(\Omega)} \big) \\
 \nonumber
   & \le c \norm{\grad\u}{\H}^{3/2} \norm{\u}{\H}^{1/2} 
   \le \frac18 \norm{\grad\u}{\H}^2 + c \norm{\u}{\H}^2,
\end{align}
having used regularity \eqref{rego:u:s} both for $\u_1$ and for $\u_2$, Sobolev's embeddings,
and the Poincar\'e-Wirtinger inequality (recall that both $\u_1$ and $\u_2$ have zero spatial
mean). 

\smallskip

Next, we multiply \eqref{eq: A-C pi diff} by $-\lap\varphi$ to deduce
\begin{align}\label{co:52}
 \frac12\derivt\norm{\grad\varphi}{H} 
  + \epsilon \norm{\lap\varphi}{H}^2 
  + \sigma \normH{\grad\vp}^2
 & = - (\pi,\lap\fhi) 
 + \big(\u_1 \cdot \grad\varphi + \u \cdot \grad\varphi_2,\lap\fhi\big).
\end{align}
Correspondingly, we multiply \eqref{eq: A-C pi_t diff} by $\pi$. Standard
manipulations give
\begin{equation}\label{co:53}
 \frac{\ka}{2}\derivt\norm{\pi}{H} 
  + \norm{\pi}{H}^2 
   = - \delta \big( \u_1 \cdot \nabla \pi + \u \cdot \nabla \pi_2 , \pi \big) 
   + (\pi,\lap\fhi) 
   - \big(f(\varphi_1)- f(\varphi_2),\pi\big).
\end{equation}
Combining the previous two relations we obtain
\begin{align}\label{co:54}
 & \frac12\derivt\norm{\grad\varphi}{H} 
  + \frac\kappa2\derivt\norm{\pi}{H} 
  + \epsilon \norm{\lap\varphi}{H}^2 
  + \sigma \normH{\grad\vp}^2
  + \norm{\pi}{H}^2 \\
 \nonumber
  & \mbox{}~~~~~
   = \big(\u_1 \cdot \grad\varphi ,\lap\fhi\big)
   + \big(\u \cdot \grad\varphi_2,\lap\fhi\big)
   - \delta \big( \u_1 \cdot \nabla \pi, \pi \big) 
   - \delta \big( \u \cdot \nabla \pi_2 , \pi \big) 
    - \big(f(\varphi_1)- f(\varphi_2),\pi\big)\\
 \nonumber 
  & \mbox{}~~~~~
   = : J_1 + J_2 + J_3 + J_4 + J_5,
\end{align}
and we need to control the quantities $J_i$ on the \rhs. First of all,
\begin{align}\label{st:J1}
  | J_1 | & \le \norm{\u_1}{L^6(\Omega)} \norm{\grad\varphi}{L^3(\Omega)} \norm{\lap\varphi}{H} \\
 \nonumber
  & \le \norm{\u_1}{\V} \norm{\grad\varphi}{H}^{1/2} \norm{\lap\varphi}{H}^{3/2} 
    \le \frac{\epsilon}6 \norm{\lap\varphi}{H}^2 + c_\epsilon \norm{\grad\varphi}{H}^2,
\end{align}
where we used interpolation, regularity \eqref{rego:u:s} for $\u_1$, and elliptic
regularity results for $\varphi$. Next, using the fact that $\u$ has zero spatial 
mean and the regularity \eqref{rego:fhi:s} for 
$\fhi_2$, we obtain
\begin{align}\label{st:J2}
  | J_2 | & \le \norm{\u}{L^3(\Omega)} \norm{\grad\varphi_2}{L^6(\Omega)} \norm{\lap\varphi}{H} \\
 \nonumber
  & \le \norm{\u}{\H}^{1/2}\norm{\grad\u}{\H}^{1/2} \norm{\varphi_2}{H^2(\Omega)} \norm{\lap\varphi}{H} 
    \le \frac18 \norm{\grad\u}{\H}^2 + \frac{\epsilon}6 \norm{\lap\varphi}{H}^2 + c_\epsilon \norm{\u}{\H}^2.
\end{align}
On the other hand, by incompressibility, 
\begin{equation}\label{st:J3}
  J_3 = - \frac\delta2 \ints{ \u_1 \cdot \nabla \pi^2 } = 0,
\end{equation}
whereas the subsequent term, for $\delta>0$, can be controlled only in case we have the additional 
regularity \eqref{rego:add} (here for $i=2$):
\begin{align}\label{st:J4}
  | J_4 | & \le \delta \norm{\u}{L^6(\Omega)} \norm{\grad\pi_2}{L^3(\Omega)} \norm{\pi}{H} \\
 \nonumber
  & \le \frac18\norm{\grad\u}{\H}^2 + c \delta^2 \norm{\grad\pi_2}{L^3(\Omega)}^2 \norm{\pi}{H}^2.
\end{align}
Finally, using \eqref{powerlike}, we have
\begin{equation}\label{st:J5}
  | J_5 |  \le c \ints{ ( 1 + |\varphi_2|^p + |\varphi_1|^p ) |\pi| |\fhi| }
   \le c \normH{\fhi}^2 + c \normH{\pi}^2.
\end{equation}

\smallskip

Now, we can take the sum of \eqref{co:51} and \eqref{co:54}. Using
\eqref{st:K1}-\eqref{st:K2} and \eqref{st:J1}-\eqref{st:J5},
we then deduce
\begin{align}\label{co:54b}
 & \frac12 \derivt\norm{\u}{\H}^2 
  + \frac12 \norm{\grad\u}{\H}^2
  + \frac12\derivt\norm{\grad\varphi}{H} 
  + \frac{\ka}{2} \derivt\norm{\pi}{H} 
  + \frac\epsilon2 \norm{\lap\varphi}{H}^2 
  + \sigma \normH{\grad\vp}^2 \\
  %  + \norm{\pi}{H}^2  
 \nonumber
  & \mbox{}~~~~~
   \le c_\epsilon \norm{\grad\fhi}{H}^2
   + c_\epsilon \norm{\u}{\H}^2
   + c \big( 1 + \delta^2 \norm{\grad\pi_2}{L^3(\Omega)}^2 \big) \norm{\pi}{H}^2
   + c \normH{\fhi}^2.
\end{align}
To control the last term on the \rhs, we need to test \eqref{eq: A-C pi diff} with $\varphi$.
This procedure yields
\begin{align}\label{eq: stima FASE}
  & \frac12\derivt\norm{\varphi}{H}^2 
   + \epsilon \norm{\grad\varphi}{H}^2
   + \sigma \normH{\vp}^2\\
  \nonumber 
  & \mbox{} = \ints{\pi\varphi} 
    - \ints{\big (\u_1\cdot\grad\varphi ) \fhi}
    - \ints { (\u \cdot\grad\varphi_2) \varphi} =: H_1 + H_2 + H_3.
\end{align}
Now, $H_2$ is readily seen to be $0$ thanks to incompressibility. On the other hand, 
it is easy to check that
\begin{align}\label{st:H1}
  | H_1 | & \le c \normH{\vp}^2 + c \normH{\pi}^2,\\
 \label{st:H2}
  | H_3 | & \le \norm{\u}{L^4(\Omega)} \norm{\grad\fhi_2}{L^4(\Omega)} \norm{\varphi}{H}
    \le \frac14 \norm{\grad\u}{\H}^2 + c \normH{\vp}^2.
\end{align}
Hence, summing \eqref{eq: stima FASE} to \eqref{co:54b},
neglecting some positive quantities on the \lhs, and taking
\eqref{st:H1}-\eqref{st:H2} into account, we finally arrive at
\begin{align}\label{co:54c}
 & \frac12 \derivt\norm{\u}{\H}^2 
  + \frac14 \norm{\grad\u}{\H}^2
  + \frac12\derivt\norm{\varphi}{V}^2
  + \frac{\ka}{2} \derivt\norm{\pi}{H}^2
  + \frac\epsilon2 \norm{\lap\varphi}{H}^2 
  + \sigma \normH{\grad\vp}^2 \\
 \nonumber
  & \mbox{}~~~~~
   \le c_\epsilon \norm{\grad\fhi}{H}^2
   + c_\epsilon \norm{\u}{\H}^2
   + c \big( 1 + \delta^2 \norm{\grad\pi_2}{L^3(\Omega)}^2 \big) \norm{\pi}{H}^2
   + c \normH{\fhi}^2.
\end{align}
Hence, exploiting in the case $\delta>0$ the additional regularity assumption
\eqref{rego:add}, we can use Gronwall's lemma in the above relation
to deduce that $(\u_1, \varphi_1, \pi_1)$ coincides with $(\u_2, \varphi_2, \pi_2)$
over the whole of $(0,T)$, which actually concludes the proof of the theorem.
\begin{remark}\label{rem:dipcon}{\rm
 It is clear that, in the case when the initial data for
 $(\u_1, \varphi_1, \pi_1)$ and $(\u_2, \varphi_2, \pi_2)$ do not coincide
 with each other, then one can obtain from \eqref{co:54c} a continuous dependence
 estimates in the norms specified in \eqref{hp:init2}. Namely,
 one has 
 \begin{align}\label{dip:con}
  & \norm{\u_1(t) - \u_2(t)}{\H}^2 
    + \norm{\varphi_1(t) - \fhi_2(t)}{V}^2 
    + \norm{\pi_1(t) - \pi_2(t)}{H}^2\\
 \nonumber   
  & \mbox{}~~~~~
   \le C(T) \Big( \norm{\u_1(0) - \u_2(0)}{\H}^2 
    + \norm{\varphi_1(0) - \fhi_2(0)}{V}^2 
    + \norm{\pi_1(0) - \pi_2(0)}{H}^2 \Big)
 \end{align}
 for any $t\in[0,T]$, with the constant $C(T)$ on the \rhs\ depending 
 on the ``strong'' norms of the two
 solutions specified in \eqref{rego:u:s}-\eqref{rego:pi:s}
 (and in \eqref{rego:add} in the case $\delta>0$). }
\end{remark}

%%%%%%%%%%%%%%%%%%%%%%%%%%%%%%%%%%%%%%%%%%%%%%%%%%%%%%%%%%%%%%%%%%%%%%%%%%%%%%%%%%%%%%%%%%%%%%%%%%%%%%%%%%%%%%%%%%%%%%%%%%%%%%

\section{Galerkin approximation}
\label{sec:Gal app}

In this section we present a possible construction of a sequence $(\u_n,\fhi_n,\pi_n)$
of approximate solutions by means of a Faedo-Galerkin scheme. Since the procedure is rather
standard and follows an approach already used for similar models (see, e.g., 
\cite[p.~284]{Miranville-Caginalp} or \cite{ColliN-S--A-C}),
we will only give some highlights leaving the remaining details to the reader.

That said, we let $\Hz$ and $\Vz$ denote the subspaces, respectively, of 
$\H$ and of $\V$, consisting of the functions having zero spatial mean.
Then we can define the Stokes operator as an unbounded linear operator on $\Hz$ 
by setting $A = -P\lap$, $D(A)=H^2(\Omega)\cap \Vz$, where $P:L^2(\Omega)\to \H$ is the Leray 
projector \cite{temam1984navier}. Notice that 
\begin{equation*}
  (A\u,\v)=(\!(\u,\v)\!) 
   =(\grad\u,\grad\v)\quad \perogni \u\in D(A), \quad \perogni \v \in \Vz.
\end{equation*}
Moreover, the operator $A$ is positive and self-adjoint on $\Hz$. Hence, we can take
as a Galerkin base of $\Hz$ the family $\{\w_i\}_{i\in\NN}$ of the (linearly independent
and properly normalized) eigenfunctions of $A$. Next, noting as $I$ the identity operator of $H$, we consider the
unbounded linear operator $B$ of $H$ defined as $B:= I - \Delta$ with domain $D(B) = H^2(\Omega)$ 
(note that $\Omega$-periodicity is still implicitly assumed also at this level). 
Hence, $B$ is also positive and self-adjoint and we can take as a Galerkin base of $H$
the family $\{\tau_i\}_{i\in\NN}$ of the (normalized and linearly independent)
eigenfunctions of $B$. For any $m\in \NN$, 
we can then define the $m$-dimensional subspaces $W_m\ede \spann\{\w_1, \dots, \w_m\}\subset \Hz$ 
and $T_m \ede \spann\{\tau_1, \dots, \tau_m\}\subset H$. We also denote as $P_m: \Hz \to W_m$ 
and $\Pi_m: H \to T_m$ the orthogonal projectors onto $W_m$ and $T_m$, respectively.

Finally, it is convenient to replace $f$ 
with a suitable regularization $f_n$ depending on a further approximation parameter $n$.
More precisely, we assume that $f_n\in C^1(\RR)$ with 
\begin{equation}\label{hp:fn}
  |f_n(r)|+|f_n'(r)| \le C_n \quad \perogni n\in\NN,~r\in\RR,
\end{equation}
with the constants $C_n>0$ of course going to infinity as $n\nearrow \infty$. An explicit
expression of $f_n$ can be easily constructed by suitably 
truncating $f$ outside some bounded interval $I_n$ 
increasing with respect to $n$, for example $I_n=[-n,n]$. Then we obtain as a byproduct
that $f_n$ converges to $f$ uniformly on compact sets of $\RR$ as $n\nearrow \infty$.

With this machinery at hand, we can look for a Faedo-Galerkin solution 
of the form
\begin{equation*}
  \u_m = \sum_{j=1}^{m}u^m_j(t)\w_j, \qquad  
   \vp_m = \sum_{j=1}^{m}\vp^m_j(t) \tau_j, \qquad 
   \pi_m = \sum_{j=1}^{m}\pi^m_j(t) \tau_j,
\end{equation*}
satisfying the following discretized system:
\begin{align}\label{NS-gal}
  & (\partial_t \u_m,\w)
   + (\u_m\cdot\grad\u_m,\w)
   + (A \u_m,\w)
   = (\grad\vp_m \otimes \grad\vp_m,\grad\w),
    \quad \perogni \w\in V_m, \\
 \label{pi-gal}
   & \pi_m = \partial_t\vp_m + \Pi_m (\u_m\cdot\grad\vp_m) 
    + \epsilon B \vp_m + ( \sigma - \epsilon) \vp_m,\\
 \label{fhi-gal}
   & \ka \partial_t \pi_m + \pi_m
     + \delta \Pi_m (\u_m\cdot \grad\pi_m) 
     + B \fhi_m - \fhi_m
     + \Pi_m (f_n(\vp_m)) = 0,
\end{align}
complemented with the initial conditions    
\begin{equation}\label{init-gal}
  \u_m|_{t=0} = P_m \u_0, \qquad
  \fhi_m|_{t=0} = \Pi_m \fhi_0, \qquad
  \pi_m|_{t=0} = \Pi_m \pi_0.
\end{equation}
Since all nonlinear terms in the above system
have at least a locally Lipschitz dependence on their arguments, it 
turns out that existence of a local in time solution to
\eqref{NS-gal}-\eqref{fhi-gal} with the initial conditions
\eqref{init-gal} is a consequence of the classical
Cauchy theorem for ODE's. It is also worth stressing that the above constructed
approximate solution, despite being identified by the sole subscript
$m$, depends in fact on both approximation parameters $n$ and $m$. 

For any fixed $n\in\NN$, we will now let $m$ tend to infinity so to
obtain a solution $(\u_n,\vp_n,\pi_n)$ depending only on $n$. 
To this aim, we just reproduce the a priori bounds obtained in the previous section
working now on the Faedo-Galerkin scheme. Indeed, all the computations
can be repeated just with some small technical difference mainly related
to the presence of a regularized nonlinear function $f_n$ in place of 
the original $f$. The only point we would like to remark refers to
the use of \eqref{pi-gal} to compute the product term now taking the form
$(\Pi_m (f_n(\vp_m)),\pi_m)$ (compare with \eqref{giu:12}).
Indeed, in the present setting this gives
rise to the quantity $(\Pi_m (f_n(\vp_m)),\Pi_m (\u_m\cdot\grad\vp_m))$
which does not necessarily vanish because {\sl both factors}\/ are 
projected onto the finite dimensonal subspace $T_m$. On the other hand,
we can easily see that
\begin{align}\label{gal-11}
  \big( \Pi_m (f_n(\vp_m)),\Pi_m (\u_m\cdot\grad\vp_m) \big)
   & = \big( \Pi_m (f_n(\vp_m)), \u_m\cdot\grad\vp_m \big)\\
 \nonumber  
   & = \big( \Pi_m (f_n(\vp_m)) - f_n(\vp_m) , \u_m\cdot\grad\vp_m \big)\\
 \nonumber  
   & \le c_n \| \u_m \|_{\H} \| \grad \vp_m \|_H,
\end{align}
where we also used incompressibility and condition \eqref{hp:fn}.
Note that the additional contribution on the \rhs, at {\sl fixed}\/
$n\in\NN$, can be controlled using Gronwall (giving rise to 
an additional quantity in the approximate energy inequality).
On the other hand, as we take $m\nearrow \infty$, the projection 
operators disappear and we do no longer face any additional term when 
repeating the estimates on $(\u_n,\vp_n,\pi_n)$.

It is also worth observing that, as we take $m\nearrow \infty$, we can
rely on a set of a priori estimates that are uniform with respect to the time
variable. Hence, as a consequence of standard extension arguments
it will turn out that, despite Galerkin solutions $(\u_m,\vp_m,\pi_m)$ 
might be defined only on small time intervals, the solutions $(\u_n,\vp_n,\pi_n)$
obtained in the limit can be thought to be defined for every $t\in(0,\infty)$. 
As a consequence, we can take $(\u_n,\vp_n,\pi_n)$ as a regularized
solution to our system depending on the sole parameter $n$. Moreover,
the a priori estimates obtained before are also uniform
with respect to $n\in\NN$. Hence, to let $n\nearrow \infty$, we can proceed as 
described in Subsection~\ref{subsec:wss},
with the sole differences related to the occurrence of the 
regularized nonlinearity $f_n$. Nevertheless, even here, the needed
modifications are almost straightforward. Indeed, using the 
strong (hence a.e.~pointwise) convergence of $\fhi_n$ (cf.~\eqref{conv:22})
and the fact that $f_n\to f$ uniformly on compact subsets of $\RR$, 
it is easy to realize that, in place of \eqref{conv:26}, there holds
\begin{equation} \label{conv:26n}     
  f_n(\vp_n) \to f(\vp) \quext{weakly star in }\, L^\infty(0,T;L^{\frac{p+2}{p+1}}(\Omega))
   \quext{and strongly in }\, L^1(0,T;L^{1}(\Omega)).
\end{equation}
Then, the rest of the argument works up to minor adaptations.

%
%
%
%
%
%%%%%%%%%%%%%%%%%%%%%%%%%%%%%%%%%%%%%%%%%%%%%%%%%%%%%%%%%%%%%%%%%%%%%%%%%%%%%%%%%%%%%%%%%%%%%%%%%%%%%%%%%%%%%%%%
%
%
%
%
%
%
%
%
%
%
%
%
%

\section{Conclusive remarks}
\label{sec:conrem}

In this part, we provide some additional considerations on the studied model
and in particular we discuss some of the choices on coefficients and regularizing terms
we operated for the purpose of obtaining a tractable mathematical problem.

%%%%%%%%%%%%%%%%%%%%%%%%%%%%%%%%%%%%%%%%%%%%%%%%%%%%%%%%%%%%%%%%%%%%%%%%%%%%%%%%%%%%%%%%%%%%%%%%%%%%%%%%%%%%%%%%%%%%%%%%%%%%%%

\subsection{The model without regularizing terms}

Our system \eqref{eq: N-S}-\eqref{eq: A-C pi_t} can be seen
as a regularized version of a more basic model which can be 
written as follows:
\begin{align} \label{NS11}
  & \u_t + \u \cdot \grad \u + \grad p - \lap \u = - \dive(\grad\vp\otimes\grad\vp),\\
 \label{incom11}
  & \dive \u =0, \\ 
 \label{pi11}
  & \pi = \pi_0 := \vp_t + \u \cdot \grad\vp,\\
 \label{fhi11}
  & \kappa( \pi_t + \u \cdot \grad\pi) + \pi - \lap\vp + f(\vp) = 0,
\end{align}
The above system corresponds to \eqref{eq: N-S}-\eqref{eq: A-C pi_t} with the 
choices $\sigma=\epsilon=0$ and $\kappa = \delta$. It is worth noting that, neglecting
the role of the velocity, i.e., taking $\u  = 0$, \eqref{pi11}-\eqref{fhi11}
reduce to the damped wave equation (or, equivalently, to the hyperbolic
relaxation of the Allen-Cahn equation with ``small'' relaxation coefficient $\kappa$).
When the contribution of $\u$ is included, proceeding as in Subsec.~\ref{subsec:en}
one can easily obtain that (hypothetical) solutions to 
\eqref{NS11}-\eqref{fhi11} satisfy the energy law 
\begin{gather}\label{ene11}
  \ddt \calE + \calD_0 = 0, \quad\text{where} \\
 \label{ene12}
  \calE = \frac12 \normHd{\u}^2
     + \frac\kappa2 \normH{\pi}^2 
     + \frac12 \| \nabla \fhi \|_H^2
     + \ints{ F(\fhi) }, \qquad 
  \calD_0 = \normHd{\grad\u}^2
     + \normH{\pi}^2.
\end{gather}
Physically speaking the above relation states that the derivative of 
the energy $\calE$ is given by minus the dissipation terms $\calD_0$.
Moreover, the four summands in the energy can be respectively
interpreted as the kinetic-macroscopic, kinetic-microscopic,
interface, and configuration energies. On the other hand, the 
two summands in $\calD_0$ represent the dissipation due to 
viscosity effects in the flow and to the process of phase transition,
respectively. Relation 
\eqref{ene11} represents the First Principle of Thermodynamics,
stating that, in absence of external sources, a part of the 
energy $\calE$ is gradually converted into heat.
In particular, one obtains the finiteness of the {\sl dissipation
integrals}:
\begin{equation} \label{diss22}
   \int_0^{+\infty} \calD_0(t) \,\dit < + \infty.
\end{equation}
\begin{remark}[nonisothermal case]{\rm
 If one wants to include temperature into the model, 
 a nonisothermal version of \eqref{NS11}-\eqref{fhi11} could be written,
 following the lines, e.g., of \cite{ERS1}, in the form
 \begin{align} \label{NS12}
   & \u_t + \u \cdot \grad \u + \grad p - \lap \u = - \dive(\grad\vp\otimes\grad\vp),\\
  \label{incom12}
   & \dive \u =0, \\ 
  \label{pi12}
   & \pi = \pi_0 := \vp_t + \u \cdot \grad\vp,\\
  \label{fhi12}
   & \kappa( \pi_t + \u \cdot \grad\pi) + \pi - \lap\vp + f(\vp) = \vt,\\
  \label{vt12}
   & \vt_t + \u \cdot \nabla \vt + \pi \vt - \dive( a(\vt) \nabla \vt ) 
    = D_0 := | \nabla \u |^2 + | \pi |^2,
 \end{align}
 where $\vt>0$ is the {\sl absolute}\/ temperature
 and $a(\vt)> 0$ represents the (possibly constant) heat conductivity.
 Of course other nonisothermal effects may be also included in \eqref{NS12}-\eqref{vt12},
 like for instance a temperature dependent viscosity in \eqref{NS12}. In particular
 the model written above is just the simplest one having some physical significance.
 In this setting, if we repeat the procedure leading to the energy inequality
 but add now also the integral of \eqref{vt12} over $\Omega$, noting
 that the integral of $D_0$ equals $\calD_0$, we obtain
 (still for periodic boundary conditions) the modified energy balance
 \begin{equation} \label{ene13}
    \ddt \Big( \calE + \io \vt \, \dix \Big)
     = 0, 
 \end{equation}
 where now the energy is conserved in view of the fact that also thermal energy is 
 now considered. In addition to that, testing \eqref{vt12} by $-\vt^{-1}$
 one obtains the {\sl entropy production inequality}\/ corresponding to the 
 validity of the Second Principle of Thermodynamics.
 }
\end{remark}
It is finally worth observing that relation \eqref{ene11} (or its integral
version \eqref{diss22}) prescribes that some amount of energy is converted into heat, 
but in itself it does not suffice to quantify the amount of energy that is released.
In other words it does not imply the {\sl uniform}\/ dissipativity property
leading to existence of an absorbing set. Actually, to achieve this,
we need to mimick the procedure in the proof in Subsec.~\ref{subsec:diss}
exploiting the fact that we have taken here $\delta=\kappa$. Actually, 
testing \eqref{pi11} by $\kappa \pi$, \eqref{fhi11} by $\fhi$, and using
incompressibility, we deduce
\begin{equation} \label{diss11}
  \kappa \ddt \io \pi \fhi\, \dix 
   + \io \pi \fhi \, \dix
   + ( - \Delta \fhi + f(\fhi) , \fhi )
   = \kappa \| \fhi \|^2.
\end{equation}
Then, repeating the procedure given in Subsec.~\ref{subsec:diss}, i.e.,
summing (a positive multiple of) \eqref{diss11} to \eqref{ene11}, one obtains
uniform dissipativity. This fact suggests that, from the energetic viewpoint, 
the case $\delta=\kappa$
is physically more sound. In other words, the correct relaxation of the 
Allen-Cahn equation with transport effect should be obtained by adding
the second {\sl material}\/ derivative of $\fhi$, as done in \eqref{fhi11}
(and not, for instance, just its second partial derivative $\fhi_{tt}$). 
Of course, this choice has a mathematical drawback, namely, estimating the term 
$\vu \cdot \nabla \pi$ may be complicated, especially in the setting
of strong solutions. This is exactly the issue which led us to consider also,
in our mathematical results, the case when $\kappa>0$ and $\delta=0$.

%%%%%%%%%%%%%%%%%%%%%%%%%%%%%%%%%%%%%%%%%%%%%%%%%%%%%%%%%%%%%%%%%%%%%%%%%%%%%%%%%%%%%%%%%%%%%%%%%%%%%%%%%%%%%%%%%%%%%%%%%%%%%%

\subsection{Introduction of regularizing terms}

In the previous subsection we have seen that system \eqref{NS11}-\eqref{fhi11}
complies with the basic laws of Thermodynamics. Moreover, at least in absence
of external source terms, solution trajectories satisfy a uniform dissipation principle. 
In this sense, the above model can be seen as a physically reasonable and very
natural extension of the usual models for two-phase fluids (like those considered, e.g.,
in \cite{Ab,AbelsRoger'09,Boyer,ColliN-S--A-C}). On the other hand,
as explained in the introduction, system \eqref{NS11}-\eqref{fhi11} 
seems impervious to a rigorous mathematical analysis,
even in two dimensions of space. In other words, if one wants to prove
(at least) existence of weak solutions, some modification is mandatory.
Indeed, the simplest and more natural strategy to obtain a mathematically 
tractable PDE system consists in the addition of regularizing terms, 
and this is what we did for the sake of proving our results. Of course,
regularizing the model may affect the underlying physical aspects,
and now we would like to discuss a bit this issue.

Since the main mathematical difficulties arise in connection with the regularity 
of $\fhi$, it appears to be natural to add regularizing terms either in \eqref{pi11}
or in \eqref{fhi11}, and we actually made the first choice. Generally speaking,
this corresponds to replacing \eqref{pi11} with
\begin{equation} \label{pi13}
  \pi = \pi_0 + \pi_{\reg}, \quext{where, still, }\, \pi_0 = \vp_t + \u \cdot \grad\vp,
\end{equation}
and now $\pi_{\reg}$ corresponds to the regularization term(s). The main reason for adding
regularization terms in \eqref{pi11} (rather than, for instance, in \eqref{fhi11})
stands in the fact that this choice affects 
in a minimal way the other equations and, in particular, it does not generate 
(artificial) additional forces on the \rhs\ of \eqref{NS11}. Indeed, whatever is the
choice of $\pi_{\reg}$, the strategy leading to the energy estimate (or, in other
words, the variational structure of the system) remains essentially the same.
Indeed, a simple check shows that \eqref{ene11} is now replaced by 
\begin{equation} \label{ene14}
  \ddt \calE + \calD_0 + \calD_{\reg} = 0,
    \quext{where }\, \calD_{\reg} = \big( \pi_{\reg} , - \Delta \fhi + f(\fhi) \big)
\end{equation}
and one needs to choose $\pi_{\reg}$ in such a way that the additional
dissipation $\calD_{\reg}$ has reasonable properties. In our setting, we actually took
\begin{equation} \label{pireg}
   \pi_{\reg} = - \epsilon \Delta \fhi + \sigma \fhi,
\end{equation}
in such a way that
\begin{equation} \label{Dreg}
  \calD_{\reg} = \epsilon \| \Delta \fhi \|^2 
    + \io \big( \sigma + \epsilon f'(\fhi) \big) | \nabla \fhi |^2 \, \dix
    + \sigma \io f(\fhi) \fhi \,\dix.
\end{equation}
Here the main issue stands in the nonconvexity of $F$ (or, equivalently, on the 
nonmonotonicity of $f$) implying that, particularly 
when $\sigma = 0$, $\calD_{\reg}$ may 
be unbounded from below, (possibly) providing a growth, rather than a dissipation of 
energy in that case, i.e.~a clearly nonphysical behavior (note, however, that 
in the case when $\kappa=\delta$, this issue may be somehow overcome as in 
Subsec.~\ref{subsec:diss} or in \eqref{diss11}).

On the other hand, if $\epsilon$ and $\sigma$ are {\sl both}\/ strictly
positive, even though the quantity $\calD_0 + \calD_{\reg}$ may be nonpositive
(due to nonmonotonicity of $f$), it is however true that, up to a positive
constant, $\calD_0 + \calD_{\reg}$ essentially controls the energy from above
(namely, the above argument implies the inequality \eqref{eq: 04b4}).
Hence, \eqref{Dreg} {\sl in itself}\/ provides a uniform dissipativity principle,
with no need of adding \eqref{diss11}. For this reason, such a choice works also
in the case when $\kappa\not=\delta$, i.e., \eqref{fhi11} is replaced by
\begin{equation} \label{fhidiff}
  \kappa \pi_t + \delta \u \cdot \grad\pi + \pi - \lap\vp + f(\vp) = 0.
\end{equation}
It is nevertheless true that, with the position \eqref{pireg},
one loses the finiteness of the dissipation integrals \eqref{diss22}
and, for this reason, in the case \eqref{pireg} the compatibility of the model
with the Energy conservation principle cannot be considered to be complete.
It is worth noting that, taking instead
\begin{equation} \label{pireg2}
  \pi_{\reg} = \epsilon \big( - \Delta \fhi + f(\fhi) \big),
\end{equation}
we would have 
\begin{equation} \label{Dreg2}
  \calD_{\reg} = \epsilon \| - \Delta \fhi + f(\fhi) \|^2 
\end{equation}
whence the occurrence of global dissipation integrals would be kept in 
this case:
\begin{equation} \label{diss23}
   \int_0^{+\infty} \big( \calD_0(t) + \calD_{\reg(t)} \big) \,\dit < + \infty
\end{equation}
As noted in the introduction, the model with the choice \eqref{pireg2} 
may be studied in a forthcoming work.

%%%%%%%%%%%%%%%%%%%%%%%%%%%%%%%%%%%%%%%%%%%%%%%%%%%%%%%%%%%%%%%%%%%%%%%%%%%%%%%%%%%%%%%%%%%%%%%%%%%%%%%%%%%%%%%%%%%%%%%%%%%%%%

\subsection{Neglecting the second material derivative}

A second simplification we operated in order to prove our mathematical results
(particularly those devoted to strong solutions)
is given by the choice $\delta = 0$, that corresponds to replacing 
\eqref{fhi11} with
\begin{equation} \label{fhi00}
  \kappa \pi_t + \pi - \lap\vp + f(\vp) = 0.
\end{equation}
As noted above, on the one hand this choice has a poor motivation from the physical point
of view (we saw that relaxation should be provided through the material derivative
of $\pi$, i.e., assuming $\kappa=\delta$). 
On the other hand, it can be justified, at least partially,
as a linearization argument which leads to neglecting some terms that 
depend at least quadratically on some ``small'' quantities. Indeed, if
expanded, the (neglected) second order transport term can be 
written as 
\begin{equation} \label{transp}
  \delta \u \cdot \nabla \pi 
   = \delta \u \cdot \nabla (\fhi_t + \u \cdot \nabla \fhi),
\end{equation}
and, at least for small flow velocity $\u$, the above quantity can be though 
to be negligible compared to other terms. This is of course an ``ad-hoc'' simplification,
but it is a procedure commonly performed in the mathematical literature 
on phase transition models, the most striking example being probably given by
non-isothermal phase-field model. In that framework most of the mathematical
works are devoted to models that are linearized around the critical temperature 
losing in that case the validity of the Second Principle of Thermodynamics. This
is the case, for instance, of the so-called ``Caginalp phase-field model''
\cite{Caginalp'86} which is studied, along with its many variants,
in a huge number of mathematical papers. This model can be obtained
from other types of (thermodynamically consistent, but mathematically 
more difficult) models exactly by neglecting some quadratic terms 
depending on ``small'' quantities, see for instance \cite{FPR}
and the references therein.

%%%%%%%%%%%%%%%%%%%%%%%%%%%%%%%%%%%%%%%%%%%%%%%%%%%%%%%%%%%%%%%%%%%%%%%%%%%%%%%%%%%%%%%%%%%%%%%%%%%%%%%%%%%%%%%%%%%%%%%%%%%%%%

\bigskip

\noindent%
{\bf Acknowledgments.}~~The present paper benefits from the support of the Italian
MIUR-PRIN Grant 2015PA5MP7 ``Calculus of Variations'' and of the GNAMPA 
(Gruppo Nazionale per l'A\-na\-li\-si Ma\-te\-ma\-ti\-ca, la Pro\-ba\-bi\-li\-t\`a e le 
loro Ap\-pli\-ca\-zio\-ni) of INdAM (I\-sti\-tu\-to Na\-zio\-na\-le di 
Alta Mate\-ma\-ti\-ca) for GS. GF has been supported
by the Austrian Science Fund (FWF) grant W1245 and by Vienna Doctoral School 
of Mathematics.

% % % % % % % %BIBLIOGRAFIA % % % % % % % % % %

\pagestyle{plain}

\end{document}